# New Results Based On Riemann Hypothesis Is Tenable

**Kaida Shi**

Department of Mathematics, Zhejiang Ocean University,

Zhoushan 316004, Zhejiang Province, China

**Abstract** Starting from the Euler's identity, the author improved Riemann's results, discovered the relationship between the Riemann $\zeta$ function $\zeta(s)$ and the prime function $\omega(s)$, and obtained two new corollaries under the Riemann hypothesis is tenable. By using these corollaries, the author found the relationship between the $p_m$, its subscript $m$ and the $t_m$, and obtained a Complete table of primes (first 1000 primes, less than 7919).

## 1. Two new corollaries under the Riemann hypothesis is true

In the well-known paper "Ueber die Anzahl der Primzahlen unter einer gegebenen Große"[1], German mathematician B. Riemann used the Euler's identity

$$\zeta(s) = \sum_{n=1}^{\infty} n^{-s} = \prod_p (1 - p^{-s})^{-1},$$

which embodies the relationship between all the natural numbers and all the primes.

Although Riemann took natural logarithm on two sides of the above identity and obtained:

$$\log \zeta(s) = \log \prod_p (1 - p^{-s})^{-1} = -\sum_p \log(1 - p^{-s})$$

he developed it and obtained

$$\begin{aligned} &\log \prod_p (1 - p^{-s})^{-1} = -\sum_p \log(1 - p^{-s}) \\ &= \sum_p p^{-s} + \frac{1}{2}\sum_p p^{-2s} + \frac{1}{3}\sum_p p^{-3s} + \cdots\cdots + \frac{1}{\mu}\sum_p p^{-\mu s} + \cdots\cdots. \end{aligned} \tag{1}$$

But, he did not pay attention to the following result:

$$\begin{aligned} &\log \zeta(s) = (\zeta(s) - 1) - \frac{(\zeta(s)-1)^2}{2} + \frac{(\zeta(s)-1)^3}{3} - \cdots\cdots + \\ &+ (-1)^{\mu+1}\frac{(\zeta(s)-1)^{\mu}}{\mu} + \cdots\cdots \\ &= -(1-\zeta(s)) - \frac{(1-\zeta(s))^2}{2} - \frac{(1-\zeta(s))^3}{3} - \cdots\cdots - \frac{(1-\zeta(s))^{\mu}}{\mu} - \cdots\cdots \end{aligned} \tag{2}$$

As L. Euler used contrast method to proved audaciously $\sum_{n=1}^{\infty}\frac{1}{n^2}=\frac{\pi^2}{6}$ (**please refer to the Appendix below**), we contrast both expressions (1) and (2), can obtain the following equation group:

$$\begin{cases}\sum_p p^{-s}=-(1-\zeta(s)),\\ \sum_p p^{-2s}=-(1-\zeta(s))^2,\\ \sum_p p^{-3s}=-(1-\zeta(s))^3,\\ \cdots\cdots\cdots\cdots\cdots\cdots\cdots\cdots\cdots\cdots\cdots,\\ \sum_p p^{-\mu s}=-(1-\zeta(s))^{\mu},\\ \cdots\cdots\cdots\cdots\cdots\cdots\cdots\cdots\cdots\cdots\cdots.\end{cases} \tag{3}$$

From above equation group, we obtain:

$$\left(\sum_p p^{-s}\right)^{\mu}=(-1)^{\mu}(-1)(-1)(1-\zeta(s))^{\mu}$$
$$=(-1)^{\mu+1}(-(1-\zeta(s))^{\mu})=(-1)^{\mu+1}\cdot\sum_p p^{-\mu s},$$

namely, when $\mu$ is odd number,

$$\left(\sum_p p^{-s}\right)^{\mu}=\sum_p p^{-\mu s},$$

when $\mu$ is even number,

$$\left(\sum_p p^{-s}\right)^{\mu}=-\sum_p p^{-\mu s}.$$

Substituting above results into the equation group (3), we obtain:

$$\begin{cases}\sum_p p^{-s}=\zeta(s)-1,\\ -\left(\sum_p p^{-s}\right)^2=-(\zeta(s)-1)^2,\\ \left(\sum_p p^{-s}\right)^3=(\zeta(s)-1)^3,\\ \cdots\cdots\cdots\cdots\cdots\cdots\cdots\cdots\cdots\cdots\cdots,\\ (-1)^{\mu+1}\left(\sum_p p^{-s}\right)^{\mu}=(-1)^{\mu+1}(\zeta(s)-1)^{\mu}\\ \cdots\cdots\cdots\cdots\cdots\cdots\cdots\cdots\cdots\cdots\cdots.\end{cases}$$

So, the solution of above equation group is

$$\sum_p p^{-s} = \zeta(s) - 1. \tag{4}$$

If we call the $\sum_p p^{-s}$ as the prime function, then, (4) shows the relationship between the Riemann $\zeta$ function $\zeta(s) = \sum_{n=1}^{\infty} n^{-s}$ and the prime function $\sum_p p^{-s}$ .

**Corollary 1** If Riemann hypothesis is tenable, namely

$$\zeta(s) = \sum_{n=1}^{\infty} n^{-s} = 0, \qquad (s = \frac{1}{2} + it,\ \ t \neq 0) \tag{5}$$

Substituting (5) into (4), we obtain:

$$\sum_p p^{-(\frac{1}{2}+it)} = -1.$$

Because

$$\sum_p p^{-(\frac{1}{2}+it)} = \sum_p p^{-\frac{1}{2}} \cdot e^{-it\log p}$$

$$= \sum_p p^{-\frac{1}{2}} \cdot (\cos(t\log p) - i\sin(t\log p)) = -1 = -1 + i \cdot 0,$$

therefore, we have

$$\sum_p p^{-\frac{1}{2}} \cdot \cos(t\log p) = -1. \tag{6}$$

**Corollary 2** Considering the prime set adds the non-prime set equals the natural number set, we can obtain:

$$\begin{cases} \omega(s) = \sum_p p^{-s} = -1; \\ \zeta(s) = \sum_{n=1}^{\infty} n^{-s} = 0; \qquad (s = \frac{1}{2} + it,\ \ t \neq 0) \\ \lambda(s) = \sum_c c^{-s} = 1. \end{cases}$$

where, $\omega(s) = \sum_p p^{-s} = -1$ is called the **prime function equation**, and $\lambda(s) = \sum_c c^{-s} = 1$ is called the **non-prime function equation**.

## 2. The relationship between the primes and the non-trivial zeroes

Writing the equation (6) as the inner product of two infinite-dimensional vectors, we obtain:

$$\omega(s) = (\frac{1}{\sqrt{p_1}},\ \frac{1}{\sqrt{p_2}},\ \frac{1}{\sqrt{p_3}},\ \cdots\cdots,\ \frac{1}{\sqrt{p_m}},\ \cdots\cdots) \bullet (\cos(t_1 \log p_1),$$
$$\cos(t_2 \log p_2),\ \cos(t_3 \log p_3),\ \cdots\cdots,\ \cos(t_m \log p_m),\ \cdots\cdots) = -1. \quad (7)$$

On the other hand, we have

$$\omega(s) = (\frac{1}{\sqrt{p_1}},\ \frac{1}{\sqrt{p_2}},\ \frac{1}{\sqrt{p_3}},\ \cdots\cdots,\ \frac{1}{\sqrt{p_m}},\ \cdots\cdots) \bullet$$
$$\bullet (\frac{-\sqrt{p_1}}{1\cdot 2},\ \frac{-\sqrt{p_2}}{2\cdot 3},\ \frac{-\sqrt{p_3}}{3\cdot 4},\ \cdots\cdots,\ \frac{-\sqrt{p_m}}{m(m+1)},\ \cdots\cdots) = -1, \quad (8)$$

Considering the consistence between above two expressions (7) and (8), we can obtain:

$$\cos(t_m \log p_m) = \frac{-\sqrt{p_m}}{m(m+1)},$$

namely

$$m^2 + m + \frac{\sqrt{p_m}}{\cos(t_m \log p_m)} = 0. \quad (*)$$

This shows the relationship between the prime $p_m$, its subscript $m$ and the $t_m$ (the imagination part of the non-trivial zero of Riemann $\zeta$ function $\zeta(s)$ ).

We can use the equation (*) to draw up a complete table of primes (see the Appendix below).

## 3. The proof about the uniqueness of the solution of the equation (*)

Suppose that for a certain $m$, the equation (*) has two solutions, they are respectively $(p^1_{,m},\ t^1_m)$ and $(p^2_{,m},\ t^2_m)$. Substituting them into the equation (*), we obtain:

$$m^2 + m + \frac{\sqrt{p^1_m}}{\cos(t^1_m \log p^1_m)} = 0, \quad (9)$$

$$m^2 + m + \frac{\sqrt{p^2_m}}{\cos(t^2_m \log p^2_m)} = 0, \quad (10)$$

Taking (9) to minus (10), we have

$$\frac{\sqrt{p^1_m}}{\cos(t^1_m \log p^1_m)} - \frac{\sqrt{p^2_m}}{\cos(t^2_m \log p^2_m)} = 0,$$

namely

$$\frac{\sqrt{p_m^1}}{\cos(t_m^1 \log p_m^1)} = \frac{\sqrt{p_m^2}}{\cos(t_m^2 \log p_m^2)},$$

therefore

$$\frac{\sqrt{p_m^1}}{\sqrt{p_m^2}} = \frac{\cos(t_m^1 \log p_m^1)}{\cos(t_m^2 \log p_m^2)}.$$

Suppose that

$$\frac{\sqrt{p_m^1}}{\sqrt{p_m^2}} = \frac{\cos(t_m^1 \log p_m^1)}{\cos(t_m^2 \log p_m^2)} = l, \quad (l \neq 1),$$

we will have

$$p_m^1 = l^2 p_m^2.$$

If $p_m^1 = l^2 p_m^2$ is still an integer, then

1) If $l^2$ is an integer, then, because $p_m^1 = l^2 p_m^2$, therefore, $p_m^1$ is a composite number;

2) If $l^2$ is a rational fraction, because $p_m^1 = l^2 p_m^2$ is an integer, therefore, $p_m^2$ is a composite number;

3) If $l^2$ is an irrational number, because $p_m^1 = l^2 p_m^2$, therefore, $p_m^1$ is also an irrational number.

But in fact, $p_m^1$ and $p_m^2$ are all primes, therefore, the $l$ must equal to 1. So, we have

$p_m^1 = p_m^2$ and $t_m^1 = t_m^2$.

## 4. The change scope of the $t_m$

Solving the equation (*), we can obtain:

$$m = \frac{-1 + \sqrt{1 - \dfrac{4\sqrt{p_m}}{\cos(t_m \log p_m)}}}{2}. \tag{11}$$

Obviously, when $t_m \log p_m = \frac{\pi}{2}$, the function $\frac{\sqrt{p_m}}{\cos(t_m \log p_m)} = \infty$, therefore $t_m \log p_m = \frac{\pi}{2}$ is the interrupted point of the function $\frac{\sqrt{p_m}}{\cos(t_m \log p_m)}$. From $t_m \log p_m = \frac{\pi}{2}$, we can obtain $t_m = \frac{\pi}{2 \log p_m}$. When the $p_m$ is determined, we denote

$t_0 = \frac{\pi}{2\log p_m}$.

If we take $t_m < t_0$, so $1 - \frac{4\sqrt{p_m}}{\cos(t_m \log p_m)} < 0$, but $m$ is a positive integer, therefore, it can't be taken the value from the left side of the $t_0$. When we take $t_m > t_0$, because $\cos(t_m \log p_m) < \cos(t_0 \log p_m) = \cos\frac{\pi}{2} = 0$, therefore, $1 - \frac{4\sqrt{p_m}}{\cos(t_m \log p_m)} > 0$.

Because

$$\frac{\partial m}{\partial t_m} = \frac{-\sqrt{p_m} \cdot \sin(t_m \log p_m) \cdot \log p_m}{\sqrt{1 - \frac{4\sqrt{p_m}}{\cos(t_m \log p_m)}} \cdot \cos^2(t_m \log p_m)},$$

therefore, when $\frac{\partial m}{\partial t_m} = 0$, we have

$$\sin(t_m \log p_m) = 0.$$

So, the arrest point of the function $m(t_m)$ is

$$t_m \log p_m = \pi,$$

we denote the arrest point as

$$t_m^0 = \frac{\pi}{\log p_m}.$$

Obviously, we can know

$$t_0 = \frac{\pi}{2\log p_m} < t_m < t_m^0 = \frac{\pi}{\log p_m} = 2t_0,$$

namely,

$$t_m \in (t_0,\ 2t_0).$$

On the other hand, from (*), we have

$$t_m = \frac{\cos^{-1}\left(\frac{-\sqrt{p_m}}{m(m+1)}\right)}{\log p_m}.$$

Comparing $t_m$ and $t_0$, we obtain

$$\lim_{m\to\infty}\frac{t_m}{t_0}=\lim_{m\to\infty}\frac{2\cos^{-1}\left(\frac{-\sqrt{p_m}}{m(m+1)}\right)}{\pi}=1,$$

$$\lim_{m\to\infty}(t_m-t_0)=\lim_{m\to\infty}\frac{1}{\log p_m}\cdot\lim_{m\to\infty}\left(\cos^{-1}\left(\frac{-\sqrt{p_m}}{m(m+1)}\right)-\frac{\pi}{2}\right)=0.$$

We can know

$$t_m\in(t_0,\ 2t_0).$$

# 5. The recurrence formula for solving $p_{m+k}$ or $p_m$

Suppose that for a pair $m_1$ and $m_2$, we have a pair $(p_{m+k},\ t_{m+k})$ and the $(p_m,\ t_m)$. According to the expression (11), we obtain

$$m_1-m_2=\frac{\sqrt{1-\frac{4\sqrt{p_{m+k}}}{\cos(t_{m+k}\log p_{m+k})}}}{2}-\frac{\sqrt{1-\frac{4\sqrt{p_m}}{\cos(t_m\log p_m)}}}{2}=k,$$

where, $k$ is the distance (in subscript) between two primes $p_{m+k}$ and $p_m$.

After putting above expression in order, we can obtain:

$$\frac{\sqrt{p_{m+k}}}{\cos(t_{m+k}\log p_{m+k})}=\frac{\sqrt{p_m}}{\cos(t_m\log p_m)}-k\cdot\sqrt{1-\frac{4\sqrt{p_m}}{\cos(t_m\log p_m)}}-k^2;$$

or

$$\frac{\sqrt{p_m}}{\cos(t_m\log p_m)}=\frac{\sqrt{p_{m+k}}}{\cos(t_{m+k}\log p_{m+k})}+k\cdot\sqrt{1-\frac{4\sqrt{p_{m+k}}}{\cos(t_{m+k}\log p_{m+k})}}+k^2.$$

So long as the three values $p_m$, $t_m$ and $k$ are determined, we can obtain the $p_{m+k}$ and $t_{m+k}$. Conversely, so long as the three values $p_{m+k}$, $t_{m+k}$ and $k$ are determined, we can obtain the $p_m$ and $t_m$.

# 6. Conclusion

We improved Riemann's results, discovered the relationship between the Riemann $\zeta$ function $\zeta(s)$ and the prime function $\omega(s)$, and obtained two new corollaries based on Riemann hypothesis is tenable. By using these corollaries, we found the relationship between the $p_m$, its subscript $m$ and the $t_m$, and drew up a complete table of primes (first 1000 primes, less than 7927), but because of the equation (*) is an indefinite equation, therefore, when either $p_m$ or $m$ is fixed, if we can obtain the corresponded another by regulating the $t_m$ in its change scope, then in fact, the expression (*) is just the **term formula** for solving the prime $p_m$. On the other hand, if we can find the largest prime $p_m$ less than $x$, by using the expression (*), we can obtain $m$ (the number of primes less than $x$). So, the prime distribution problem will be solved.

**Acknowledgements** I am grateful to all those who expressed their interest and support to the results of this research. Especially, I thank Director, Dr. Bernd Siebert at the Mathematisches Institut der Georg-August-Universität, Göttingen, German for his help during the supplement of the duplication of B. Riemann's original text "Ueber die Anzahl der Primzahlen unter einer gegebenen Große" and Associate Professor Tianhui Liu at Zhejiang Ocean Institute, China for his help during the translation of B. Riemann's original text (from German translate into Chinese).

## References

1. Riemann, B., Ueber die Anzahl der Primzahlen unter einer gegebenen Große, Ges. Math. Werke und Wissenschaftlicher Nachlaß , 2, Aufl, 1859, pp 145~155.

2. Lehman R. S., On the difference $\pi(x)$-li$x$, *Acta Arith.*, 11(1966). pp.397~410.

3. Hardy, G. H., Littlewood, J. E., Some problems of "patitio numerorum" III: On the expression of a number as a sum of primes, *Acta. Math.*, 44 (1923). pp.1~70.

4. Hardy, G. H., Ramanujan, S., Asymptotic formula in combinatory analysis, *Proc. London Math. Soc.*, (2) 17 (1918). pp. 75~115.

5. E. C. Titchmarsh, *The Theory of the Riemann Zeta Function*, Oxford University Press, New York, 1951.

6. Morris Kline, *Mathematical Thought from Anoient to Modern Times*，Oxford University Press, New York, 1972.

7. A. Selberg, *The zeta and the Riemann Hypothesis*, Skandinaviske Mathematiker Kongres, 10 (1946).

8. Lou S. T., Wu D. H., *Riemann hypothesis*, Shenyang: Liaoning Education Press, 1987. pp.152~154.

9. Pan C. D., Pan C. B., *Goldbach Conjection*, Beijing: Science Press, 1981. pp.1~18; pp.119~147.

10. Hua L. G., *A Guide to the Number Theory*, Beijing: Science Press, 1957. Pp234~263.

# Appendix

## How was L. Euler proved $\sum_{n=1}^{\infty}\frac{1}{n^2}=\frac{\pi^2}{6}$ by using the contrast method ?

In the end of 17 century, Swiss mathematician J. Bernoulli wanted to know what the infinite series

$$1+\frac{1}{2^2}+\frac{1}{3^2}+\frac{1}{4^2}+\frac{1}{5^2}+\cdots\cdots$$

equals to ? If somebody solved the problem and announces him, he would appreciate him.

First, L. Euler (1707--1783) supposed that the equation

$$a_0-a_1x^2+a_2x^4-\cdots+(-1)^n a_n x^{2n}=0 \qquad (1)$$

has $2n$ different roots:

$$\alpha_1,\ -\alpha_1,\ \alpha_2,\ -\alpha_2,\ \cdots,\ \alpha_n,\ -\alpha_n\,.$$

According to the factorization of the polynomial, he obtained:

$$\begin{aligned}&a_0-a_1x^2+a_2x^4-\cdots+(-1)^n a_n x^{2n}\\&=a_0\left(1-\frac{x^2}{\alpha_1^2}\right)\left(1-\frac{x^2}{\alpha_2^2}\right)\cdots\left(1-\frac{x^2}{\alpha_n^2}\right)\\&=a_0-a_0\left(\frac{1}{\alpha_1^2}+\frac{1}{\alpha_2^2}+\cdots+\frac{1}{\alpha_n^2}\right)x^2+\cdots.\end{aligned}$$

Obviously,

$$a_1=a_0\left(\frac{1}{\alpha_1^2}+\frac{1}{\alpha_2^2}+\cdots+\frac{1}{\alpha_n^2}\right). \qquad (2)$$

Secondary, considering the equation:

$$\sin x=0\,,$$

it equals to:

$$x-\frac{x^3}{3!}+\frac{x^5}{5!}-\frac{x^7}{7!}+\cdots=0\,. \qquad (3)$$

L.Euler considered the left hand of above equation has infinite terms and infinite-th power, therefore, it must has infinite roots:

$$0,\ \pi,\ -\pi,\ 2\pi,\ -2\pi,\ 3\pi,\ -3\pi,\ \cdots.$$

Removing the root 0 and use $x$ to divide two sides of (3), he obtained:

$$1-\frac{x^2}{3!}+\frac{x^4}{5!}-\frac{x^6}{7!}+\cdots=0 \qquad (4)$$

The roots of above equation (4) are respectively:

$$\alpha_1=\pi,\ \ \alpha_1^{'}=-\pi,\ \ \alpha_2=2\pi,\ \ \alpha_2^{'}=-2\pi,\ \ \alpha_3=3\pi,\ \ \alpha_3^{'}=-3\pi,\ \ \cdots.$$

At the moment, L. Euler took (1) and (4) to **contrast audaciously**, and obtained:

$$\begin{aligned}&a_0-a_1x^2+a_2x^4-\cdots+(-1)^n a_n x^{2n}\\&=a_0\left(1-\frac{x^2}{\alpha_1^2}\right)\left(1-\frac{x^2}{\alpha_2^2}\right)\cdots\left(1-\frac{x^2}{\alpha_n^2}\right),\\&\frac{\sin x}{x}=1-\frac{x^2}{3!}+\frac{x^4}{5!}-\frac{x^6}{7!}+\cdots\\&=1\cdot\left(1-\frac{x^2}{(\pm\pi)^2}\right)\left(1-\frac{x^2}{(\pm2\pi)^2}\right)\left(1-\frac{x^2}{(\pm3\pi)^2}\right)\cdots,\end{aligned}$$

where

$$a_0=1,\ \ a_1=\frac{1}{3!},\ \ a_2=\frac{1}{5!},\ \ \cdots\cdots.$$

Also, according to (2), he obtained:

$$\frac{1}{3!}=\frac{1}{\pi^2}+\frac{1}{4\pi^2}+\frac{1}{9\pi^2}+\cdots..$$

Taking $\pi^2$ to mulitplies two sides of above expression, he obtained:

$$\frac{\pi^2}{6}=1+\frac{1}{4}+\frac{1}{9}+\cdots,$$

namely

$$\frac{\pi^2}{6}=1+\frac{1}{2^2}+\frac{1}{3^2}+\cdots.$$

# Explanation:

1) In following Table, the Colour Blue denotes the numerlcal value of imagine part $t_m$ of non-trival zero of Riemann $\zeta$ Function $\zeta(s)$ which is corresponded by the prime $p_m$; the Colour Orange denotes the numerlcal value of imagine part $t_\eta$ of non-trival zero of Riemann $\zeta$ Function $\zeta(s)$ which is corresponded by the composite number $c_\eta$.

2) The formulas for calculating the $t_m$ and $t_\eta$ are respectively:

$$t_m = \frac{\cos^{-1}\left(\frac{-\sqrt{p_m}}{m(m+1)}\right)}{\log p_m}$$

and

$$t_\eta = \frac{\cos^{-1}\left(\frac{\sqrt{c_\eta}}{\eta(\eta+1)}\right)}{\log c_\eta}.$$

# The Complete Table of Primes and Composite Numbers

| $N$ | $m$ | $\eta$ | Blue $t_m$ Orange $t_\eta$ | $N$ | $m$ | $\eta$ | Blue $t_m$ Orange $t_\eta$ |
|---|---|---|---|---|---|---|---|
| 1 | | 1 | $\infty$ | 50 | | 35 | 0.4000958932 |
| 2 | 1 | | 3.3992701060 | 51 | | 36 | 0.3981445299 |
| 3 | 2 | | 1.6963574120 | 52 | | 37 | 0.3962467539 |
| 4 | | 2 | 0.8879495235 | 53 | 16 | | 0.4023796406 |
| 5 | 3 | | 1.0924502700 | 54 | | 38 | 0.3925405083 |
| 6 | | 3 | 0.7619479172 | 55 | | 39 | 0.3907941424 |
| 7 | 4 | | 0.8754118490 | 56 | | 40 | 0.3890922917 |
| 8 | | 4 | 0.6871552980 | 57 | | 41 | 0.3874331213 |
| 9 | | 5 | 0.6693120589 | 58 | | 42 | 0.3858148930 |
| 10 | | 6 | 0.6494581617 | 59 | 17 | | 0.3876633904 |
| 11 | 5 | | 0.7012720514 | 60 | | 43 | 0.3826503189 |
| 12 | | 7 | 0.6072251650 | 61 | 18 | | 0.3876633904 |
| 13 | 6 | | 0.6459186318 | 62 | | 44 | 0.3796386060 |
| 14 | | 8 | 0.5755105912 | 63 | | 45 | 0.3782068406 |
| 15 | | 9 | 0.5641512425 | 64 | | 46 | 0.3768069478 |
| 16 | | 10 | 0.5534267161 | 65 | | 47 | 0.3754377642 |
| 17 | 7 | | 0.5804327797 | 66 | | 48 | 0.3740981819 |
| 18 | | 11 | 0.5323362189 | 67 | 19 | | 0.3787050658 |
| 19 | 8 | | 0.5540524507 | 68 | | 49 | 0.3714723682 |
| 20 | | 12 | 0.5147739249 | 69 | | 50 | 0.3702171449 |
| 21 | | 13 | 0.5076706478 | 70 | | 51 | 0.3689874663 |
| 22 | | 14 | 0.5009504656 | 71 | 20 | | 0.3732065227 |
| 23 | 9 | | 0.5179753118 | 72 | | 52 | 0.3665746726 |
| 24 | | 15 | 0.4878402246 | 73 | 21 | | 0.3704244073 |
| 25 | | 16 | 0.4822841958 | 74 | | 53 | 0.3642581105 |
| 26 | | 17 | 0.4770061522 | 75 | | 54 | 0.3631464340 |
| 27 | | 18 | 0.4719900824 | 76 | | 55 | 0.3620555102 |
| 28 | | 19 | 0.4672194834 | 77 | | 56 | 0.3609846948 |
| 29 | 10 | | 0.4810304852 | 78 | | 57 | 0.3599333699 |
| 30 | | 20 | 0.4580019237 | 79 | 22 | | 0.3635156856 |
| 31 | 11 | | 0.4697131298 | 80 | | 58 | 0.3578669641 |
| 32 | | 21 | 0.4497029754 | 81 | | 59 | 0.3568715659 |
| 33 | | 22 | 0.4460002345 | 82 | | 60 | 0.3558933727 |
| 34 | | 23 | 0.4424484680 | 83 | 23 | | 0.3592122065 |
| 35 | | 24 | 0.4390388809 | 84 | | 61 | 0.3539692645 |
| 36 | | 25 | 0.4357631293 | 85 | | 62 | 0.3530405336 |
| 37 | 12 | | 0.4458141449 | 86 | | 63 | 0.3521270830 |
| 38 | | 26 | 0.4294097593 | 87 | | 64 | 0.3512284903 |
| 39 | | 27 | 0.4265072638 | 88 | | 65 | 0.3503443489 |
| 40 | | 28 | 0.4237079054 | 89 | 24 | | 0.3534526181 |
| 41 | 13 | | 0.4324638276 | 90 | | 66 | 0.3486038452 |
| 42 | | 29 | 0.4182678755 | 91 | | 67 | 0.3477613334 |
| 43 | 14 | | 0.4259351409 | 92 | | 68 | 0.3469317593 |
| 44 | | 30 | 0.4132096674 | 93 | | 69 | 0.3461147875 |
| 45 | | 31 | 0.4108675190 | 94 | | 70 | 0.3453100941 |
| 46 | | 32 | 0.4085975832 | 95 | | 71 | 0.3445173669 |
| 47 | 15 | | 0.4154036691 | 96 | | 72 | 0.3437363042 |
| 48 | | 33 | 0.4041695107 | 97 | 25 | | 0.3466774192 |
| 49 | | 34 | 0.4021033399 | 98 | | 73 | 0.3421973595 |

| $N$ | $m$ | $\eta$ | Blue $t_m$ Orange $t_\eta$ | $N$ | $m$ | $\eta$ | Blue $t_m$ Orange $t_\eta$ |
|---|---|---|---|---|---|---|---|
| 99 | | 74 | 0.3414499745 | 152 | | 116 | 0.3124851235 |
| 100 | | 75 | 0.3407131282 | 153 | | 117 | 0.3120802629 |
| 101 | 26 | | 0.3434607769 | 154 | | 118 | 0.3116790442 |
| 102 | | 76 | 0.3392604818 | 155 | | 119 | 0.3112814107 |
| 103 | 27 | | 0.3418152806 | 156 | | 120 | 0.3108873069 |
| 104 | | 77 | 0.3378480375 | 157 | | | 0.3124271024 |
| 105 | | 78 | 0.3371608870 | 158 | | 121 | 0.3101067376 |
| 106 | | 79 | 0.3364828486 | 159 | | 122 | 0.3097229629 |
| 107 | 28 | | 0.3388815946 | 160 | | 123 | 0.3093425058 |
| 108 | | 80 | 0.3351449370 | 161 | 37 | 124 | 0.3089653167 |
| 109 | 29 | | 0.3373863987 | 162 | | 125 | 0.3085913467 |
| 110 | | 81 | 0.3338418945 | 163 | | | 0.3100684523 |
| 111 | | 82 | 0.3332069791 | 164 | | 126 | 0.3078504228 |
| 112 | | 83 | 0.3325800089 | 165 | | 127 | 0.3074858810 |
| 113 | 30 | | 0.3346936830 | 166 | | 128 | 0.3071243721 |
| 114 | | 84 | 0.3313418883 | 167 | 38 | | 0.3085350152 |
| 115 | | 85 | 0.3307379914 | 168 | | 129 | 0.3064079756 |
| 116 | | 86 | 0.3301413804 | 169 | | 130 | 0.3060553541 |
| 117 | | 87 | 0.3295519051 | 170 | | 131 | 0.3057055942 |
| 118 | | 88 | 0.3289694197 | 171 | 39 | 132 | 0.3053586555 |
| 119 | | 89 | 0.3283937822 | 172 | | 133 | 0.3050144983 |
| 120 | | 90 | 0.3278248550 | 173 | | | 0.3063704960 |
| 121 | | 91 | 0.3272625041 | 174 | | 134 | 0.3043322960 |
| 122 | | 92 | 0.3267065990 | 175 | | 135 | 0.3039962961 |
| 123 | | 93 | 0.3261570133 | 176 | | 136 | 0.3036629261 |
| 124 | | 94 | 0.3256136235 | 177 | 40 | 137 | 0.3033321501 |
| 125 | | 95 | 0.3250763097 | 178 | | 138 | 0.3030039331 |
| 126 | | 96 | 0.3245449554 | 179 | | | 0.3043085611 |
| 127 | 31 | | 0.3266093677 | 180 | | 139 | 0.3023531558 |
| 128 | | 97 | 0.3234947179 | 181 | | | 0.3035965569 |
| 129 | | 98 | 0.3229807096 | 182 | | 140 | 0.3017123151 |
| 130 | | 99 | 0.3224722188 | 183 | 41 | 141 | 0.3013964595 |
| 131 | 32 | | 0.3244248382 | 184 | | 142 | 0.3010829639 |
| 132 | | 100 | 0.3214668164 | 185 | 42 | 143 | 0.3007717974 |
| 133 | | 101 | 0.3209744001 | 186 | | 144 | 0.3004629300 |
| 134 | | 102 | 0.3204871044 | 187 | | 145 | 0.3001563321 |
| 135 | | 103 | 0.3200048378 | 188 | | 146 | 0.2998519746 |
| 136 | | 104 | 0.3195275107 | 189 | | 147 | 0.2995498289 |
| 137 | 33 | | 0.3213891550 | 190 | | 148 | 0.2992498671 |
| 138 | | 105 | 0.3185833246 | 191 | | | 0.3004605526 |
| 139 | 34 | | 0.3203388907 | 192 | | 149 | 0.2986548237 |
| 140 | | 106 | 0.3176581744 | 193 | | | 0.2998110688 |
| 141 | | 107 | 0.3172044742 | 194 | | 150 | 0.2980682736 |
| 142 | | 108 | 0.3167552144 | 195 | 43 | 151 | 0.2977788825 |
| 143 | | 109 | 0.3163103201 | 196 | | 152 | 0.2974915160 |
| 144 | | 110 | 0.3158697180 | 197 | 44 | | 0.2986023209 |
| 145 | | 111 | 0.3154333366 | 198 | | 153 | 0.2969213008 |
| 146 | | 112 | 0.3150011061 | 199 | | | 0.2979842046 |
| 147 | | 113 | 0.3145729581 | 200 | | 154 | 0.2963589671 |
| 148 | | 114 | 0.3141488260 | 201 | 45 | 155 | 0.2960814105 |
| 149 | 35 | | 0.3158475580 | 202 | | 156 | 0.2958057351 |
| 150 | | 115 | 0.3133092183 | 203 | 46 | 157 | 0.2955319187 |
| 151 | 36 | | 0.3149160280 | 204 | | 158 | 0.2952599394 |

| $N$ | $m$ | $\eta$ | Blue $t_m$ Orange $t_\eta$ | $N$ | $m$ | $\eta$ | Blue $t_m$ Orange $t_\eta$ |
|---|---|---|---|---|---|---|---|
| 205 | | 159 | 0.2949897756 | 258 | | 203 | 0.2828056706 |
| 206 | | 160 | 0.2947214062 | 259 | | 204 | 0.2826093379 |
| 207 | | 161 | 0.2944548104 | 260 | | 205 | 0.2824140285 |
| 208 | | 162 | 0.2941899677 | 261 | | 206 | 0.2822258799 |
| 209 | | 163 | 0.2939268579 | 262 | | 207 | 0.2820264420 |
| 210 | | 164 | 0.2936654612 | 263 | 56 | | 0.2828128834 |
| 211 | 47 | | 0.2947079482 | 264 | | 208 | 0.2816421991 |
| 212 | | 165 | 0.2931465410 | 265 | | 209 | 0.2814518759 |
| 213 | | 166 | 0.2928901874 | 266 | | 210 | 0.2812625213 |
| 214 | | 167 | 0.2926354703 | 267 | | 211 | 0.2810741267 |
| 215 | | 168 | 0.2923823714 | 268 | | 212 | 0.2808866834 |
| 216 | | 169 | 0.2921308729 | 269 | 57 | | 0.2816512393 |
| 217 | | 170 | 0.2918809568 | 270 | | 213 | 0.2805140185 |
| 218 | | 171 | 0.2916326059 | 271 | 58 | | 0.2812519790 |
| 219 | | 172 | 0.2913858030 | 272 | | 214 | 0.2801450810 |
| 220 | | 173 | 0.2911405312 | 273 | | 215 | 0.2799622848 |
| 221 | | 174 | 0.2908967738 | 274 | | 216 | 0.2797803907 |
| 222 | | 175 | 0.2906545145 | 275 | | 217 | 0.2795993907 |
| 223 | 48 | | 0.2916766064 | 276 | | 218 | 0.2794192771 |
| 224 | | 176 | 0.2901734284 | 277 | 59 | | 0.2801374363 |
| 225 | | 177 | 0.2899355829 | 278 | | 219 | 0.2790611222 |
| 226 | | 178 | 0.2896991721 | 279 | | 220 | 0.2788836293 |
| 227 | 49 | | 0.2906839621 | 280 | | 221 | 0.2787069926 |
| 228 | | 179 | 0.2892296404 | 281 | 60 | | 0.2794035680 |
| 229 | 50 | | 0.2901750901 | 282 | | 222 | 0.2783557224 |
| 230 | | 180 | 0.2887657103 | 283 | 61 | | 0.2790291789 |
| 231 | | 181 | 0.2885362780 | 284 | | 223 | 0.2780078100 |
| 232 | | 182 | 0.2883081923 | 285 | | 224 | 0.2778353602 |
| 233 | 51 | | 0.2892205137 | 286 | | 225 | 0.2776637241 |
| 234 | | 183 | 0.2878551047 | 287 | | 226 | 0.2774928950 |
| 235 | | 184 | 0.2876309893 | 288 | | 227 | 0.2773228661 |
| 236 | | 185 | 0.2874081660 | 289 | | 228 | 0.2771536310 |
| 237 | | 186 | 0.2871866216 | 290 | | 229 | 0.2769851829 |
| 238 | | 187 | 0.2869663432 | 291 | | 230 | 0.2768175155 |
| 239 | 52 | | 0.2878510528 | 292 | | 231 | 0.2766506225 |
| 240 | | 188 | 0.2865286966 | 293 | 62 | | 0.2773117551 |
| 241 | 53 | | 0.2873799363 | 294 | | 232 | 0.2763186572 |
| 242 | | 189 | 0.2860959991 | 295 | | 233 | 0.2761540551 |
| 243 | | 190 | 0.2858818876 | 296 | | 234 | 0.2759902025 |
| 244 | | 191 | 0.2856689682 | 297 | | 235 | 0.2758270933 |
| 245 | | 192 | 0.2854572293 | 298 | | 236 | 0.2756647216 |
| 246 | | 193 | 0.2852466592 | 299 | | 237 | 0.2755030815 |
| 247 | | 194 | 0.2850372468 | 300 | | 238 | 0.2753421672 |
| 248 | | 195 | 0.2848289809 | 301 | | 239 | 0.2751819730 |
| 249 | | 196 | 0.2846218503 | 302 | | 240 | 0.2750224931 |
| 250 | | 197 | 0.2844158443 | 303 | | 241 | 0.2748637219 |
| 251 | 54 | | 0.2852491381 | 304 | | 242 | 0.2747056539 |
| 252 | | 198 | 0.2840064347 | 305 | | 243 | 0.2745482837 |
| 253 | | 199 | 0.2836362281 | 306 | | 244 | 0.2743916057 |
| 254 | | 200 | 0.2836021448 | 307 | 63 | | 0.2750451905 |
| 255 | | 201 | 0.2834016132 | 308 | | 245 | 0.2740798935 |
| 256 | | 202 | 0.2832021438 | 309 | | 246 | 0.2739252648 |
| 257 | 55 | | 0.2840114764 | 310 | | 247 | 0.2737713072 |

| $N$ | $m$ | $\eta$ | Blue $t_m$ Orange $t_\eta$ | $N$ | $m$ | $\eta$ | Blue $t_m$ Orange $t_\eta$ |
|---|---|---|---|---|---|---|---|
| 311 | 64 | 248 | 0.2744063405 | 364 | | 292 | 0.2663273444 |
| 312 | | | 0.2734649861 | 365 | | 293 | 0.2662037063 |
| 313 | 65 | 249 | 0.2740801622 | 366 | | 294 | 0.2660805196 |
| 314 | | 250 | 0.2731616009 | 367 | 73 | | 0.2665954600 |
| 315 | | 251 | 0.2730106311 | 368 | | 295 | 0.2658352385 |
| 316 | | | 0.2728606021 | 369 | | 296 | 0.2657133910 |
| 317 | 66 | 252 | 0.2734382910 | 370 | | 297 | 0.2655919835 |
| 318 | | 253 | 0.2725620656 | 371 | | 298 | 0.2654710130 |
| 319 | | 254 | 0.2724139348 | 372 | | 299 | 0.2653504769 |
| 320 | | 255 | 0.2722664259 | 373 | 74 | | 0.2658541489 |
| 321 | | 256 | 0.2721195344 | 374 | | 300 | 0.2651104570 |
| 322 | | 257 | 0.2719732555 | 375 | | 301 | 0.2649912095 |
| 323 | | 258 | 0.2718275849 | 376 | | 302 | 0.2648723853 |
| 324 | | 259 | 0.2716825180 | 377 | | 303 | 0.2647539818 |
| 325 | | 260 | 0.2715380505 | 378 | | 304 | 0.2646359965 |
| 326 | | 261 | 0.2713941779 | 379 | 75 | | 0.2651287847 |
| 327 | | 262 | 0.2712508959 | 380 | | 305 | 0.2644010404 |
| 328 | | 263 | 0.2711082003 | 381 | | 306 | 0.2642842957 |
| 329 | | 264 | 0.2709660868 | 382 | | 307 | 0.2641679587 |
| 330 | | | 0.2708245512 | 383 | 76 | | 0.2646488391 |
| 331 | 67 | 265 | 0.2714163208 | 384 | | 308 | 0.2639362746 |
| 332 | | 266 | 0.2705428636 | 385 | | 309 | 0.2638211477 |
| 333 | | 267 | 0.2704030404 | 386 | | 310 | 0.2637064185 |
| 334 | | 268 | 0.2702637787 | 387 | | 311 | 0.2635920847 |
| 335 | | 269 | 0.2701250745 | 388 | | 312 | 0.2634781437 |
| 336 | | | 0.2699869238 | 389 | 77 | | 0.2639489036 |
| 337 | 68 | 270 | 0.2705646773 | 390 | | 313 | 0.2632512171 |
| 338 | | 271 | 0.2697119504 | 391 | | 314 | 0.2631384424 |
| 339 | | 272 | 0.2695754404 | 392 | | 315 | 0.2630260511 |
| 340 | | 273 | 0.2694394685 | 393 | | 316 | 0.2629140409 |
| 341 | | 274 | 0.2693040310 | 394 | | 317 | 0.2628024096 |
| 342 | | 275 | 0.2691691241 | 395 | | 318 | 0.2626911548 |
| 343 | | 276 | 0.2690347443 | 396 | | 319 | 0.2657185218 |
| 344 | | 277 | 0.2689008878 | 397 | 78 | | 0.2630425465 |
| 345 | | 278 | 0.2687675511 | 398 | | 320 | 0.2623594257 |
| 346 | | | 0.2686347305 | 399 | | 321 | 0.2622496565 |
| 347 | 69 | 279 | 0.2692025338 | 400 | | 322 | 0.2621402527 |
| 348 | | | 0.2683703334 | 401 | 79 | | 0.2761535776 |
| 349 | 70 | 280 | 0.2689215808 | 402 | | 323 | 0.2619223362 |
| 350 | | 281 | 0.2681079672 | 403 | | 324 | 0.2618140176 |
| 351 | | 282 | 0.2679776807 | 404 | | 325 | 0.2617060558 |
| 352 | | | 0.2678478896 | 405 | | 326 | 0.2615984489 |
| 353 | 71 | 283 | 0.2683849381 | 406 | | 327 | 0.2614911947 |
| 354 | | 284 | 0.2675895005 | 407 | | 328 | 0.2613842913 |
| 355 | | 285 | 0.2674611784 | 408 | | 329 | 0.2612777365 |
| 356 | | 286 | 0.2673333386 | 409 | 80 | | 0.2617212887 |
| 357 | | 287 | 0.2672059778 | 410 | | 330 | 0.2610654793 |
| 358 | | | 0.2670790929 | 411 | | 331 | 0.2609599600 |
| 359 | 72 | 288 | 0.2676041052 | 412 | | 332 | 0.2608547814 |
| 360 | | 289 | 0.2668264708 | 413 | | 333 | 0.2607499415 |
| 361 | | 290 | 0.2667009971 | 414 | | 334 | 0.2606454384 |
| 362 | | 377 | 0.2665759869 | 415 | | 335 | 0.2605412702 |
| 363 | | 378 | 0.2664669144 | 416 | | 336 | 0.2604374350 |

| $N$ | $m$ | $\eta$ | Blue $t_m$ Orange $t_\eta$ | $N$ | $m$ | $\eta$ | Blue $t_m$ Orange $t_\eta$ |
|---|---|---|---|---|---|---|---|
| 417 | | 337 | 0.2603339308 | 470 | | 379 | 0.2552761307 |
| 418 | | 338 | 0.2602307559 | 471 | | 380 | 0.2551880811 |
| 419 | 81 | | 0.2606677385 | 472 | | 381 | 0.2551002784 |
| 420 | | 339 | 0.2600252135 | 473 | | 382 | 0.2550127213 |
| 421 | 82 | | 0.2604512155 | 474 | | 383 | 0.2549254085 |
| 422 | | 340 | 0.2598209698 | 475 | | 384 | 0.2548383390 |
| 423 | | 341 | 0.2597194161 | 476 | | 385 | 0.2547515114 |
| 424 | | 342 | 0.2596181805 | 477 | | 386 | 0.2546649246 |
| 425 | | 343 | 0.2595172614 | 478 | | 387 | 0.2545785774 |
| 426 | | 344 | 0.2594166571 | 479 | 92 | | 0.2549304335 |
| 427 | | 345 | 0.2593163657 | 480 | | 388 | 0.2544064765 |
| 428 | | 346 | 0.2592163856 | 481 | | 389 | 0.2543208419 |
| 429 | | 347 | 0.2591167150 | 482 | | 390 | 0.2542354422 |
| 430 | | 348 | 0.2590173522 | 483 | | 391 | 0.2541502762 |
| 431 | 83 | | 0.2594371894 | 484 | | 392 | 0.2540653429 |
| 432 | | 349 | 0.2588193836 | 485 | | 393 | 0.2539806411 |
| 433 | 84 | | 0.2592289073 | 486 | | 394 | 0.2538961696 |
| 434 | | 350 | 0.2586226299 | 487 | 93 | | 0.2542426555 |
| 435 | | 351 | 0.2585247836 | 488 | | 395 | 0.2537277985 |
| 436 | | 352 | 0.2584272352 | 489 | | 396 | 0.2536440123 |
| 437 | | 353 | 0.2583299831 | 490 | | 397 | 0.2535604521 |
| 438 | | 354 | 0.2582330256 | 491 | 94 | | 0.2539000843 |
| 439 | 85 | | 0.2586346845 | 492 | | 398 | 0.2533938925 |
| 440 | | 355 | 0.2580398356 | 493 | | 399 | 0.2533110044 |
| 441 | | 356 | 0.2579437538 | 494 | | 400 | 0.2532283379 |
| 442 | | 357 | 0.2578479604 | 495 | | 401 | 0.2531458921 |
| 443 | 86 | | 0.2582409792 | 496 | | 402 | 0.2530636658 |
| 444 | | 358 | 0.2576570831 | 497 | | 403 | 0.2529816580 |
| 445 | | 359 | 0.2575621467 | 498 | | 404 | 0.2528998677 |
| 446 | | 360 | 0.2574674927 | 499 | 95 | | 0.2532344206 |
| 447 | | 361 | 0.2573731194 | 500 | | 405 | 0.2527368279 |
| 448 | | 362 | 0.2572790254 | 501 | | 406 | 0.2526556847 |
| 449 | 87 | | 0.2576646712 | 502 | | 407 | 0.2525747548 |
| 450 | | 363 | 0.2570915256 | 503 | 96 | | 0.2529028203 |
| 451 | | 364 | 0.2569982619 | 504 | | 408 | 0.2524134258 |
| 452 | | 365 | 0.2569052716 | 505 | | 409 | 0.2523331308 |
| 453 | | 366 | 0.2568125534 | 506 | | 410 | 0.2522530453 |
| 454 | | 367 | 0.2567201059 | 507 | | 411 | 0.2521731683 |
| 455 | | 368 | 0.2566279276 | 508 | | 412 | 0.2520934987 |
| 456 | | 369 | 0.2565360172 | 509 | 97 | | 0.2524160123 |
| 457 | 88 | | 0.2569154590 | 510 | | 413 | 0.2519346762 |
| 458 | | 370 | 0.2563528575 | 511 | | 414 | 0.2518556241 |
| 459 | | 371 | 0.2562617433 | 512 | | 415 | 0.2517767757 |
| 460 | | 372 | 0.2561708915 | 513 | | 416 | 0.2516981300 |
| 461 | 89 | | 0.2565424308 | 514 | | 417 | 0.2516196862 |
| 462 | | 373 | 0.2559898358 | 515 | | 418 | 0.2515414433 |
| 463 | 90 | | 0.2563528155 | 516 | | 419 | 0.2514634003 |
| 464 | | 374 | 0.2558098163 | 517 | | 420 | 0.2513855565 |
| 465 | | 375 | 0.2557202583 | 518 | | 421 | 0.2513079108 |
| 466 | | 376 | 0.2556309547 | 519 | | 422 | 0.2512304624 |
| 467 | 91 | | 0.2559865423 | 520 | | 423 | 0.2511532103 |
| 468 | | 377 | 0.2554529754 | 521 | 98 | | 0.2514724805 |
| 469 | | 378 | 0.2553644284 | 522 | | 424 | 0.2509991968 |

| $N$ | $m$ | $\eta$ | Blue $t_m$ Orange $t_\eta$ | $N$ | $m$ | $\eta$ | Blue $t_m$ Orange $t_\eta$ |
|---|---|---|---|---|---|---|---|
| 523 | 99 | | 0.2513117462 | 576 | | 471 | 0.2471148153 |
| 524 | | 425 | 0.2508459599 | 577 | 106 | | 0.2473974860 |
| 525 | | 426 | 0.2507696778 | 578 | | 472 | 0.2469801708 |
| 526 | | 427 | 0.2506935868 | 579 | | 473 | 0.2469131137 |
| 527 | | 428 | 0.2506176861 | 580 | | 474 | 0.2468462085 |
| 528 | | 429 | 0.2505419748 | 581 | | 475 | 0.2467794544 |
| 529 | | 430 | 0.2504664521 | 582 | | 476 | 0.2467128510 |
| 530 | | 431 | 0.2503911172 | 583 | | 477 | 0.2466463976 |
| 531 | | 432 | 0.2503159692 | 584 | | 478 | 0.2465800936 |
| 532 | | 433 | 0.2502410073 | 585 | | 479 | 0.2465139384 |
| 533 | | 434 | 0.2501662306 | 586 | | 480 | 0.2464479314 |
| 534 | | 435 | 0.2500916384 | 587 | 107 | | 0.2467273388 |
| 535 | | 436 | 0.2500172299 | 588 | | 481 | 0.2463162917 |
| 536 | | 437 | 0.2499430042 | 589 | | 482 | 0.2462506583 |
| 537 | | 438 | 0.2498689605 | 590 | | 483 | 0.2461853065 |
| 538 | | 439 | 0.2497950981 | 591 | | 484 | 0.2461200321 |
| 539 | | 440 | 0.2497214161 | 592 | | 485 | 0.2460549024 |
| 540 | | 441 | 0.2496479138 | 593 | 108 | | 0.2463299992 |
| 541 | 100 | | 0.2499593894 | 594 | | 486 | 0.2459250088 |
| 542 | | 442 | 0.2495013600 | 595 | | 487 | 0.2458603101 |
| 543 | | 443 | 0.2494283926 | 596 | | 488 | 0.2457957539 |
| 544 | | 444 | 0.2493556017 | 597 | | 489 | 0.2457313395 |
| 545 | | 445 | 0.2492829866 | 598 | | 490 | 0.2456670665 |
| 546 | | 446 | 0.2492105466 | 599 | 109 | | 0.2459379565 |
| 547 | 101 | | 0.2495169078 | 600 | | 491 | 0.2455388778 |
| 548 | | 447 | 0.2490661061 | 601 | 110 | | 0.2458046101 |
| 549 | | 448 | 0.2489941872 | 602 | | 492 | 0.2454112488 |
| 550 | | 449 | 0.2489224403 | 603 | | 493 | 0.2453476747 |
| 551 | | 450 | 0.2488508648 | 604 | | 494 | 0.2452842387 |
| 552 | | 451 | 0.2487794599 | 605 | | 495 | 0.2452209402 |
| 553 | | 452 | 0.2487082249 | 606 | | 496 | 0.2451577786 |
| 554 | | 453 | 0.2486371590 | 607 | 111 | | 0.2454195260 |
| 555 | | 454 | 0.2485662616 | 608 | | 497 | 0.2450318020 |
| 556 | | 455 | 0.2484955318 | 609 | | 498 | 0.2449690485 |
| 557 | 102 | | 0.2487981832 | 610 | | 499 | 0.2449064298 |
| 558 | | 456 | 0.2483544945 | 611 | | 500 | 0.2448439455 |
| 559 | | 457 | 0.2482842642 | 612 | | 501 | 0.2447815951 |
| 560 | | 458 | 0.2482141988 | 613 | 112 | | 0.2450394510 |
| 561 | | 459 | 0.2481442976 | 614 | | 502 | 0.2446572329 |
| 562 | | 460 | 0.2480745599 | 615 | | 503 | 0.2445952812 |
| 563 | 103 | | 0.2483723239 | 616 | | 504 | 0.2445334613 |
| 564 | | 461 | 0.2479354965 | 617 | 113 | | 0.2447870228 |
| 565 | | 462 | 0.2478662458 | 618 | | 505 | 0.2444101551 |
| 566 | | 463 | 0.2477971560 | 619 | 114 | | 0.2446590445 |
| 567 | | 464 | 0.2477282263 | 620 | | 506 | 0.2442873710 |
| 568 | | 465 | 0.2476594561 | 621 | | 507 | 0.2442262032 |
| 569 | 104 | | 0.2479524562 | 622 | | 508 | 0.2441651642 |
| 570 | | 466 | 0.2475223178 | 623 | | 509 | 0.2441042536 |
| 571 | 105 | | 0.2478094900 | 624 | | 510 | 0.2440434708 |
| 572 | | 467 | 0.2473858108 | 625 | | 511 | 0.2439828155 |
| 573 | | 468 | 0.2473178291 | 626 | | 512 | 0.2439222870 |
| 574 | | 469 | 0.2472500031 | 627 | | 513 | 0.2438618849 |
| 575 | | 470 | 0.2471823320 | 628 | | 514 | 0.2438016088 |

| $N$ | $m$ | $\eta$ | Blue $t_m$ Orange $t_\eta$ | $N$ | $m$ | $\eta$ | Blue $t_m$ Orange $t_\eta$ |
|---|---|---|---|---|---|---|---|
| 629 | | 515 | 0.2437414582 | 682 | | 559 | 0.2407211898 |
| 630 | | 516 | 0.2436814326 | 683 | 124 | | 0.2409382744 |
| 631 | 115 | | 0.2439281459 | 684 | | 560 | 0.2406132359 |
| 632 | | 517 | 0.2435616985 | 685 | | 561 | 0.2405594361 |
| 633 | | 518 | 0.2435020455 | 686 | | 562 | 0.2405057387 |
| 634 | | 519 | 0.2434425156 | 687 | | 563 | 0.2404521434 |
| 635 | | 520 | 0.2433831085 | 688 | | 564 | 0.2403986497 |
| 636 | | 521 | 0.2433238236 | 689 | | 565 | 0.2403452574 |
| 637 | | 522 | 0.2432646606 | 690 | | 566 | 0.2402919661 |
| 638 | | 523 | 0.2432056189 | 691 | 125 | | 0.2405065322 |
| 639 | | 524 | 0.2431466982 | 692 | | 567 | 0.2401856412 |
| 640 | | 525 | 0.2430878979 | 693 | | 568 | 0.2401326510 |
| 641 | 116 | | 0.2433319846 | 694 | | 569 | 0.2400797605 |
| 642 | | 526 | 0.2429706035 | 695 | | 570 | 0.2400269693 |
| 643 | 117 | | 0.2432103042 | 696 | | 571 | 0.2399742770 |
| 644 | | 527 | 0.2428537866 | 697 | | 572 | 0.2399216834 |
| 645 | | 528 | 0.2427955827 | 698 | | 573 | 0.2398691881 |
| 646 | | 529 | 0.2427374967 | 699 | | 574 | 0.2398167908 |
| 647 | 118 | | 0.2429733659 | 700 | | 575 | 0.2397644911 |
| 648 | | 530 | 0.2426216241 | 701 | 126 | | 0.2399769559 |
| 649 | | 531 | 0.2425638894 | 702 | | 576 | 0.2396601414 |
| 650 | | 532 | 0.2425062709 | 703 | | 577 | 0.2396081329 |
| 651 | | 533 | 0.2424487682 | 704 | | 578 | 0.2395562208 |
| 652 | | 534 | 0.2423913808 | 705 | | 579 | 0.2395044048 |
| 653 | 119 | | 0.2426239452 | 706 | | 580 | 0.2394526844 |
| 654 | | 535 | 0.2422768992 | 707 | | 581 | 0.2394010595 |
| 655 | | 536 | 0.2422198556 | 708 | | 582 | 0.2393495296 |
| 656 | | 537 | 0.2421629257 | 709 | 127 | | 0.2395595552 |
| 657 | | 538 | 0.2421061091 | 710 | | 583 | 0.2392467132 |
| 658 | | 539 | 0.2420494054 | 711 | | 584 | 0.2391954669 |
| 659 | 120 | | 0.2422787374 | 712 | | 585 | 0.2391443145 |
| 660 | | 540 | 0.2419362851 | 713 | | 586 | 0.2390932556 |
| 661 | 121 | | 0.2421616210 | 714 | | 587 | 0.2390422899 |
| 662 | | 541 | 0.2418236123 | 715 | | 588 | 0.2389914172 |
| 663 | | 542 | 0.2417674676 | 716 | | 589 | 0.2389406371 |
| 664 | | 543 | 0.2417114334 | 717 | | 590 | 0.2388899493 |
| 665 | | 544 | 0.2415666093 | 718 | | 591 | 0.2388393535 |
| 666 | | 545 | 0.2415996949 | 719 | 128 | | 0.2390473381 |
| 667 | | 546 | 0.2415439899 | 720 | | 592 | 0.2387383977 |
| 668 | | 547 | 0.2414883939 | 721 | | 593 | 0.2386880763 |
| 669 | | 548 | 0.2414329064 | 722 | | 594 | 0.2386378458 |
| 670 | | 549 | 0.2413775272 | 723 | | 595 | 0.2385877058 |
| 671 | | 550 | 0.2413222557 | 724 | | 596 | 0.2385376561 |
| 672 | | 551 | 0.2412670918 | 725 | | 597 | 0.2384876963 |
| 673 | 122 | | 0.2414905743 | 726 | | 598 | 0.2384378262 |
| 674 | | 552 | 0.2411570376 | 727 | 129 | | 0.2386434483 |
| 675 | | 553 | 0.2411021939 | 728 | | 599 | 0.2383383159 |
| 676 | | 554 | 0.2410474563 | 729 | | 600 | 0.2382887132 |
| 677 | 123 | | 0.2414905743 | 730 | | 601 | 0.2382391991 |
| 678 | | 555 | 0.2409382510 | 731 | | 602 | 0.2381897731 |
| 679 | | 556 | 0.2408838293 | 732 | | 603 | 0.2381404352 |
| 680 | | 557 | 0.2408295121 | 733 | 130 | | 0.2383433957 |
| 681 | | 558 | 0.2407752990 | 734 | | 604 | 0.2380419849 |

| $N$ | $m$ | $\eta$ | Blue $t_m$ Orange $t_\eta$ | $N$ | $m$ | $\eta$ | Blue $t_m$ Orange $t_\eta$ |
|---|---|---|---|---|---|---|---|
| 735 | | 605 | 0.2379929092 | 788 | | 650 | 0.2355094741 |
| 736 | | 606 | 0.2379439204 | 789 | | 651 | 0.2354657238 |
| 737 | | 607 | 0.2378950182 | 790 | | 652 | 0.2354200471 |
| 738 | | 608 | 0.2378462023 | 791 | | 653 | 0.2353754438 |
| 739 | 131 | | 0.2380465552 | 792 | | 654 | 0.2353309136 |
| 740 | | 609 | 0.2377487919 | 793 | | 655 | 0.2352864564 |
| 741 | | 610 | 0.2377002334 | 794 | | 656 | 0.2352420719 |
| 742 | | 611 | 0.2376517600 | 795 | | 657 | 0.2351977599 |
| 743 | 132 | | 0.2378492295 | 796 | | 658 | 0.2351535202 |
| 744 | | 612 | 0.2375550321 | 797 | 139 | | 0.2353362155 |
| 745 | | 613 | 0.2375068129 | 798 | | 659 | 0.2350652275 |
| 746 | | 614 | 0.2374586779 | 799 | | 660 | 0.2350212036 |
| 747 | | 615 | 0.2374106267 | 800 | | 661 | 0.2349772512 |
| 748 | | 616 | 0.2373626592 | 801 | | 662 | 0.2349333700 |
| 749 | | 617 | 0.2373147750 | 802 | | 663 | 0.2348895599 |
| 750 | | 618 | 0.2372669739 | 803 | | 664 | 0.2348458206 |
| 751 | 133 | | 0.2374622673 | 804 | | 665 | 0.2348021520 |
| 752 | | 619 | 0.2371715853 | 805 | | 666 | 0.2347585538 |
| 753 | | 620 | 0.2371240321 | 806 | | 667 | 0.2347150258 |
| 754 | | 621 | 0.2370765610 | 807 | | 668 | 0.2346715679 |
| 755 | | 622 | 0.2370291716 | 808 | | 669 | 0.2346281797 |
| 756 | | 623 | 0.2369818639 | 809 | 140 | | 0.2348095014 |
| 757 | 134 | | 0.2371719251 | 810 | | 670 | 0.2345415840 |
| 758 | | 624 | 0.2368874581 | 811 | 141 | | 0.2347201774 |
| 759 | | 625 | 0.2368403936 | 812 | | 671 | 0.2344552655 |
| 760 | | 626 | 0.2367934098 | 813 | | 672 | 0.2344122238 |
| 761 | 135 | | 0.2369835318 | 814 | | 673 | 0.2343692506 |
| 762 | | 627 | 0.2366996492 | 815 | | 674 | 0.2343263459 |
| 763 | | 628 | 0.2366529056 | 816 | | 675 | 0.2342835094 |
| 764 | | 629 | 0.2366062417 | 817 | | 676 | 0.2342407410 |
| 765 | | 630 | 0.2365596570 | 818 | | 677 | 0.2341980404 |
| 766 | | 631 | 0.2365131516 | 819 | | 678 | 0.2341554075 |
| 767 | | 632 | 0.2364667250 | 820 | | 679 | 0.2341128421 |
| 768 | | 633 | 0.2364203769 | 821 | 142 | | 0.2342898411 |
| 769 | 136 | | 0.2366084497 | 822 | | 680 | 0.2340278857 |
| 770 | | 634 | 0.2363278831 | 823 | 143 | | 0.2342022563 |
| 771 | | 635 | 0.2362817697 | 824 | | 681 | 0.2339431972 |
| 772 | | 636 | 0.2362357338 | 825 | | 682 | 0.2339009665 |
| 773 | 137 | | 0.2364211939 | 826 | | 683 | 0.2338588020 |
| 774 | | 637 | 0.2361438621 | 827 | 144 | | 0.2340308596 |
| 775 | | 638 | 0.2360980580 | 828 | | 684 | 0.2337746447 |
| 776 | | 639 | 0.2360523306 | 829 | 145 | | 0.2339441783 |
| 777 | | 640 | 0.2360066798 | 830 | | 685 | 0.2336907507 |
| 778 | | 641 | 0.2359611052 | 831 | | 686 | 0.2336489153 |
| 779 | | 642 | 0.2359156066 | 832 | | 687 | 0.2336071451 |
| 780 | | 643 | 0.2358701839 | 833 | | 688 | 0.2335654399 |
| 781 | | 644 | 0.2358248368 | 834 | | 689 | 0.2335237996 |
| 782 | | 645 | 0.2357795650 | 835 | | 690 | 0.2334822240 |
| 783 | | 646 | 0.2357343684 | 836 | | 691 | 0.2334407128 |
| 784 | | 647 | 0.2356892467 | 837 | | 692 | 0.2333992660 |
| 785 | | 648 | 0.2356441997 | 838 | | 693 | 0.2333578833 |
| 786 | | 649 | 0.2355992272 | 839 | 146 | | 0.2335259566 |
| 787 | 138 | | 0.2357835936 | 840 | | 694 | 0.2332752839 |

| $N$ | $m$ | $\eta$ | Blue $t_m$ Orange $t_\eta$ | $N$ | $m$ | $\eta$ | Blue $t_m$ Orange $t_\eta$ |
|---|---|---|---|---|---|---|---|
| 841 | | 695 | 0.2332340926 | 894 | | 740 | 0.2311374051 |
| 842 | | 696 | 0.2331929647 | 895 | | 741 | 0.2310994047 |
| 843 | | 697 | 0.2331519001 | 896 | | 742 | 0.2310614592 |
| 844 | | 698 | 0.2331108986 | 897 | | 743 | 0.2310235685 |
| 845 | | 699 | 0.2330699599 | 898 | | 744 | 0.2309857323 |
| 846 | | 700 | 0.2330290840 | 899 | | 745 | 0.2309479505 |
| 847 | | 701 | 0.2329882706 | 900 | | 746 | 0.2309102231 |
| 848 | | 702 | 0.2329475197 | 901 | | 747 | 0.2308725498 |
| 849 | | 703 | 0.2329068309 | 902 | | 748 | 0.2308349305 |
| 850 | | 704 | 0.2328662042 | 903 | | 749 | 0.2307973652 |
| 851 | | 705 | 0.2328256393 | 904 | | 750 | 0.2307598535 |
| 852 | | 706 | 0.2327851362 | 905 | | 751 | 0.2307223955 |
| 853 | 147 | | 0.2329522571 | 906 | | 752 | 0.2306849910 |
| 854 | | 707 | 0.2327042899 | 907 | 155 | | 0.2308383193 |
| 855 | | 708 | 0.2326639709 | 908 | | 753 | 0.2306103211 |
| 856 | | 709 | 0.2326237130 | 909 | | 754 | 0.2305730763 |
| 857 | 148 | | 0.2327886730 | 910 | | 755 | 0.2305358843 |
| 858 | | 710 | 0.2325433554 | 911 | 156 | | 0.2306873261 |
| 859 | 149 | | 0.2327059579 | 912 | | 756 | 0.2304616382 |
| 860 | | 711 | 0.2324632402 | 913 | | 757 | 0.2304246043 |
| 861 | | 712 | 0.2324232852 | 914 | | 758 | 0.2303876228 |
| 862 | | 713 | 0.2323833901 | 915 | | 759 | 0.2303506936 |
| 863 | 150 | | 0.2325439179 | 916 | | 760 | 0.2303138164 |
| 864 | | 714 | 0.2323037558 | 917 | | 761 | 0.2302769913 |
| 865 | | 715 | 0.2322640400 | 918 | | 762 | 0.2302402180 |
| 866 | | 716 | 0.2322243837 | 919 | 157 | | 0.2303902230 |
| 867 | | 717 | 0.2321847865 | 920 | | 763 | 0.2301668065 |
| 868 | | 718 | 0.2321452485 | 921 | | 764 | 0.2301301881 |
| 869 | | 719 | 0.2321057694 | 922 | | 765 | 0.2300936210 |
| 870 | | 720 | 0.2320663490 | 923 | | 766 | 0.2300571051 |
| 871 | | 721 | 0.2320269873 | 924 | | 767 | 0.2300206403 |
| 872 | | 722 | 0.2319876840 | 925 | | 768 | 0.2299842265 |
| 873 | | 723 | 0.2319484390 | 926 | | 769 | 0.2299478634 |
| 874 | | 724 | 0.2319092521 | 927 | | 770 | 0.2299115511 |
| 875 | | 725 | 0.2318701233 | 928 | | 771 | 0.2298752893 |
| 876 | | 726 | 0.2318310523 | 929 | 158 | | 0.2300240816 |
| 877 | 151 | | 0.2319906989 | 930 | | 772 | 0.2298028977 |
| 878 | | 727 | 0.2317530604 | 931 | | 773 | 0.2297667869 |
| 879 | | 728 | 0.2317141621 | 932 | | 774 | 0.2297307262 |
| 880 | | 729 | 0.2316753209 | 933 | | 775 | 0.2296947155 |
| 881 | 152 | | 0.2318329550 | 934 | | 776 | 0.2296587545 |
| 882 | | 730 | 0.2315977873 | 935 | | 777 | 0.2296228433 |
| 883 | 153 | | 0.2317532235 | 936 | | 778 | 0.2295869816 |
| 884 | | 731 | 0.2315204810 | 937 | 159 | | 0.2297343752 |
| 885 | | 732 | 0.2314819239 | 938 | | 779 | 0.2295153876 |
| 886 | | 733 | 0.2314434231 | 939 | | 780 | 0.2294796740 |
| 887 | 154 | | 0.2315969249 | 940 | | 781 | 0.2294440094 |
| 888 | | 734 | 0.2313665677 | 941 | 160 | | 0.2295896319 |
| 889 | | 735 | 0.2313282130 | 942 | | 782 | 0.2293728084 |
| 890 | | 736 | 0.2312899582 | 943 | | 783 | 0.2293372904 |
| 891 | | 737 | 0.2312517369 | 944 | | 784 | 0.2293018210 |
| 892 | | 738 | 0.2312135711 | 945 | | 785 | 0.2292264001 |
| 893 | | 739 | 0.2311754605 | 946 | | 786 | 0.2292310275 |

| $N$ | $m$ | $\eta$ | Blue $t_m$ Orange $t_\eta$ | $N$ | $m$ | $\eta$ | Blue $t_m$ Orange $t_\eta$ |
|---|---|---|---|---|---|---|---|
| 947 | 161 | | 0.2293751047 | 1000 | | 832 | 0.2273894536 |
| 948 | | 787 | 0.2291604086 | 1001 | | 833 | 0.2273565694 |
| 949 | | 788 | 0.2291251805 | 1002 | | 834 | 0.2273237275 |
| 950 | | 789 | 0.2290900003 | 1003 | | 835 | 0.2272909277 |
| 951 | | 790 | 0.2290548678 | 1004 | | 836 | 0.2272581701 |
| 952 | | 791 | 0.2290197829 | 1005 | | 837 | 0.2272254545 |
| 953 | 162 | | 0.2291623405 | 1006 | | 838 | 0.2271927808 |
| 954 | | 792 | 0.2289497377 | 1007 | | 839 | 0.2271601490 |
| 955 | | 793 | 0.2289147952 | 1008 | | 840 | 0.2271275588 |
| 956 | | 794 | 0.2288798999 | 1009 | 169 | | 0.2272613448 |
| 957 | | 795 | 0.2288450517 | 1010 | | 841 | 0.2270624880 |
| 958 | | 796 | 0.2288102504 | 1011 | | 842 | 0.2270300225 |
| 959 | | 797 | 0.2287754960 | 1012 | | 843 | 0.2269975983 |
| 960 | | 798 | 0.2287407883 | 1013 | 170 | | 0.2271298658 |
| 961 | | 799 | 0.2287061272 | 1014 | | 844 | 0.2269328585 |
| 962 | | 800 | 0.2286715126 | 1015 | | 845 | 0.2269005578 |
| 963 | | 801 | 0.2286369444 | 1016 | | 846 | 0.2268682982 |
| 964 | | 802 | 0.2286024226 | 1017 | | 847 | 0.2268360794 |
| 965 | | 803 | 0.2285679468 | 1018 | | 848 | 0.2268039013 |
| 966 | | 804 | 0.2285335172 | 1019 | 171 | | 0.2269348412 |
| 967 | 163 | | 0.2286753287 | 1020 | | 849 | 0.2267396523 |
| 968 | | 805 | 0.2284647783 | 1021 | 172 | | 0.2268689592 |
| 969 | | 806 | 0.2284304863 | 1022 | | 850 | 0.2266755653 |
| 970 | | 807 | 0.2283962398 | 1023 | | 851 | 0.2266435899 |
| 971 | 164 | | 0.2285343692 | 1024 | | 852 | 0.2266116547 |
| 972 | | 808 | 0.2283278662 | 1025 | | 853 | 0.2265797597 |
| 973 | | 809 | 0.2282937560 | 1026 | | 854 | 0.2265479047 |
| 974 | | 810 | 0.2282596910 | 1027 | | 855 | 0.2265160896 |
| 975 | | 811 | 0.2282256711 | 1028 | | 856 | 0.2264843144 |
| 976 | | 812 | 0.2281916962 | 1029 | | 857 | 0.2264525790 |
| 977 | 165 | | 0.2283303882 | 1030 | | 858 | 0.2264208833 |
| 978 | | 813 | 0.2281238639 | 1031 | 173 | | 0.2265492299 |
| 979 | | 814 | 0.2280900233 | 1032 | | 859 | 0.2263575959 |
| 980 | | 815 | 0.2280562273 | 1033 | 174 | | 0.2264843591 |
| 981 | | 816 | 0.2280224756 | 1034 | | 860 | 0.2262944664 |
| 982 | | 817 | 0.2279887683 | 1035 | | 861 | 0.2262629679 |
| 983 | 166 | | 0.2281260296 | 1036 | | 862 | 0.2262315085 |
| 984 | | 818 | 0.2279214698 | 1037 | | 863 | 0.2262000881 |
| 985 | | 819 | 0.2278878949 | 1038 | | 864 | 0.2261687068 |
| 986 | | 820 | 0.2278543640 | 1039 | 175 | | 0.2262942275 |
| 987 | | 821 | 0.2278208769 | 1040 | | 865 | 0.2261060464 |
| 988 | | 822 | 0.2277874335 | 1041 | | 866 | 0.2260747672 |
| 989 | | 823 | 0.2277540337 | 1042 | | 867 | 0.2260435552 |
| 990 | | 824 | 0.2277206774 | 1043 | | 868 | 0.2260123675 |
| 991 | 167 | | 0.2278567067 | 1044 | | 869 | 0.2259812182 |
| 992 | | 825 | 0.2276540788 | 1045 | | 870 | 0.2259501073 |
| 993 | | 826 | 0.2276208524 | 1046 | | 871 | 0.2259190348 |
| 994 | | 827 | 0.2275876692 | 1047 | | 872 | 0.2258880004 |
| 995 | | 828 | 0.2275545290 | 1048 | | 873 | 0.2258570041 |
| 996 | | 829 | 0.2275214317 | 1049 | 176 | | 0.2259816094 |
| 997 | 168 | | 0.2276560732 | 1050 | | 874 | 0.2257951117 |
| 998 | | 830 | 0.2274553494 | 1051 | 177 | | 0.2259182039 |
| 999 | | 831 | 0.2274223803 | 1052 | | 875 | 0.2257333708 |

| $N$ | $m$ | $\eta$ | Blue $t_m$ Orange $t_\eta$ | $N$ | $m$ | $\eta$ | Blue $t_m$ Orange $t_\eta$ |
|---|---|---|---|---|---|---|---|
| 1053 | | 876 | 0.2257025640 | 1106 | | 921 | 0.2241215669 |
| 1054 | | 877 | 0.2256717948 | 1107 | | 922 | 0.2240926796 |
| 1055 | | 878 | 0.2256410631 | 1108 | | 923 | 0.2240638258 |
| 1056 | | 879 | 0.2256103689 | 1109 | 186 | | 0.2241771210 |
| 1057 | | 880 | 0.2255797120 | 1110 | | 924 | 0.2240062064 |
| 1058 | | 881 | 0.2255490924 | 1111 | | 925 | 0.2239774528 |
| 1059 | | 882 | 0.2255185101 | 1112 | | 926 | 0.2239487324 |
| 1060 | | 883 | 0.2254879648 | 1113 | | 927 | 0.2239200451 |
| 1061 | 178 | | 0.2256101702 | 1114 | | 928 | 0.2238913909 |
| 1062 | | 884 | 0.2254269719 | 1115 | | 929 | 0.2238627697 |
| 1063 | 179 | | 0.2255477087 | 1116 | | 930 | 0.2238341815 |
| 1064 | | 885 | 0.2253661266 | 1117 | 187 | | 0.2239465673 |
| 1065 | | 886 | 0.2253357659 | 1118 | | 931 | 0.2237770918 |
| 1066 | | 887 | 0.2253054418 | 1119 | | 932 | 0.2237486021 |
| 1067 | | 888 | 0.2252751543 | 1120 | | 933 | 0.2237201450 |
| 1068 | | 889 | 0.2252449032 | 1121 | | 934 | 0.2236917205 |
| 1069 | 180 | | 0.2253644888 | 1122 | | 935 | 0.2236633286 |
| 1070 | | 890 | 0.2251844970 | 1123 | 188 | | 0.2237746859 |
| 1071 | | 891 | 0.2251543550 | 1124 | | 936 | 0.2336066305 |
| 1072 | | 892 | 0.2251242491 | 1125 | | 937 | 0.2235783358 |
| 1073 | | 893 | 0.2250941794 | 1126 | | 938 | 0.2235500734 |
| 1074 | | 894 | 0.2250641456 | 1127 | | 939 | 0.2235218432 |
| 1075 | | 895 | 0.2250341478 | 1128 | | 940 | 0.2234936452 |
| 1076 | | 896 | 0.2250041859 | 1129 | 189 | | 0.2236039887 |
| 1077 | | 897 | 0.2249742597 | 1130 | | 941 | 0.2234373338 |
| 1078 | | 898 | 0.2249443692 | 1131 | | 942 | 0.2234092318 |
| 1079 | | 899 | 0.2249145144 | 1132 | | 943 | 0.2233811617 |
| 1080 | | 900 | 0.2248846951 | 1133 | | 944 | 0.2233531234 |
| 1081 | | 901 | 0.2248549112 | 1134 | | 945 | 0.2233251169 |
| 1082 | | 902 | 0.2248251628 | 1135 | | 946 | 0.2232971421 |
| 1083 | | 903 | 0.2247954497 | 1136 | | 947 | 0.2232691988 |
| 1084 | | 904 | 0.2247657719 | 1137 | | 948 | 0.2232412871 |
| 1085 | | 905 | 0.2247361292 | 1138 | | 949 | 0.2232134070 |
| 1086 | | 906 | 0.2247065216 | 1139 | | 950 | 0.2231855582 |
| 1087 | 181 | | 0.2248258330 | 1140 | | 951 | 0.2231577408 |
| 1088 | | 907 | 0.2246473989 | 1141 | | 952 | 0.2231299547 |
| 1089 | | 908 | 0.2246178962 | 1142 | | 953 | 0.2231021999 |
| 1090 | | 909 | 0.2245884283 | 1143 | | 954 | 0.2230744762 |
| 1091 | 182 | | 0.2247064700 | 1144 | | 955 | 0.2230467836 |
| 1092 | | 910 | 0.2245295841 | 1145 | | 956 | 0.2230191221 |
| 1093 | 183 | | 0.2246462370 | 1146 | | 957 | 0.2229914915 |
| 1094 | | 911 | 0.2244708783 | 1147 | | 958 | 0.2229638919 |
| 1095 | | 912 | 0.2244415834 | 1148 | | 959 | 0.2229363232 |
| 1096 | | 913 | 0.2244123228 | 1149 | | 960 | 0.2229087852 |
| 1097 | 184 | | 0.2245277476 | 1150 | | 961 | 0.2228812780 |
| 1098 | | 914 | 0.2243538921 | 1151 | 190 | | 0.2229916331 |
| 1099 | | 915 | 0.2243247342 | 1152 | | 962 | 0.2228263447 |
| 1100 | | 916 | 0.2242956104 | 1153 | 191 | | 0.2229354535 |
| 1101 | | 917 | 0.2242665206 | 1154 | | 963 | 0.2227715338 |
| 1102 | | 918 | 0.2242374647 | 1155 | | 964 | 0.2227441794 |
| 1103 | 185 | | 0.2243518169 | 1156 | | 965 | 0.2227168553 |
| 1104 | | 919 | 0.2241794422 | 1157 | | 966 | 0.2226895616 |
| 1105 | | 920 | 0.2241504877 | 1158 | | 967 | 0.2226622981 |

| $N$ | $m$ | $\eta$ | Blue $t_m$ Orange $t_\eta$ | $N$ | $m$ | $\eta$ | Blue $t_m$ Orange $t_\eta$ |
|---|---|---|---|---|---|---|---|
| 1159 | | 968 | 0.2226350648 | 1212 | | 1015 | 0.2212333180 |
| 1160 | | 969 | 0.2226078616 | 1213 | 198 | | 0.2213368572 |
| 1161 | | 970 | 0.2225806885 | 1214 | | 1016 | 0.2211819594 |
| 1162 | | 971 | 0.2225535454 | 1215 | | 1017 | 0.2211563254 |
| 1163 | 192 | | 0.2226619183 | 1216 | | 1018 | 0.2211307184 |
| 1164 | | 972 | 0.2224993385 | 1217 | 199 | | 0.2212332446 |
| 1165 | | 973 | 0.2224722851 | 1218 | | 1019 | 0.2210795761 |
| 1166 | | 974 | 0.2224452615 | 1219 | | 1020 | 0.2210540498 |
| 1167 | | 975 | 0.2224182675 | 1220 | | 1021 | 0.2210285504 |
| 1168 | | 976 | 0.2223913033 | 1221 | | 1022 | 0.2210030778 |
| 1169 | | 977 | 0.2223643686 | 1222 | | 1023 | 0.2209776318 |
| 1170 | | 978 | 0.2223374634 | 1223 | 200 | | 0.2210792693 |
| 1171 | 193 | | 0.2224449867 | 1224 | | 1024 | 0.2209268107 |
| 1172 | | 979 | 0.2222837312 | 1225 | | 1025 | 0.2209014446 |
| 1173 | | 980 | 0.2222569143 | 1226 | | 1026 | 0.2208761050 |
| 1174 | | 981 | 0.2222301267 | 1227 | | 1027 | 0.2208507919 |
| 1175 | | 982 | 0.2222033684 | 1228 | | 1028 | 0.2208255051 |
| 1176 | | 983 | 0.2221766392 | 1229 | 201 | | 0.2209262660 |
| 1177 | | 984 | 0.2221499392 | 1230 | | 1029 | 0.2207750016 |
| 1178 | | 985 | 0.2221232683 | 1231 | 202 | | 0.2208746833 |
| 1179 | | 986 | 0.2220966263 | 1232 | | 1030 | 0.2207246032 |
| 1180 | | 987 | 0.2220700133 | 1233 | | 1031 | 0.2206994478 |
| 1181 | 194 | | 0.2221768161 | 1234 | | 1032 | 0.2206743185 |
| 1182 | | 988 | 0.2220168640 | 1235 | | 1033 | 0.2206492153 |
| 1183 | | 989 | 0.2219903376 | 1236 | | 1034 | 0.2206241380 |
| 1184 | | 990 | 0.2219638399 | 1237 | 203 | | 0.2207229689 |
| 1185 | | 991 | 0.2219373708 | 1238 | | 1035 | 0.2205740525 |
| 1186 | | 992 | 0.2219109305 | 1239 | | 1036 | 0.2205490531 |
| 1187 | 195 | | 0.2220167854 | 1240 | | 1037 | 0.2205240794 |
| 1188 | | 993 | 0.2218581254 | 1241 | | 1038 | 0.2204991315 |
| 1189 | | 994 | 0.2218317705 | 1242 | | 1039 | 0.2204742094 |
| 1190 | | 995 | 0.2218054441 | 1243 | | 1040 | 0.2204493129 |
| 1191 | | 996 | 0.2217791460 | 1244 | | 1041 | 0.2204244420 |
| 1192 | | 997 | 0.2217528762 | 1245 | | 1042 | 0.2203995968 |
| 1193 | 196 | | 0.2218577964 | 1246 | | 1043 | 0.2203747770 |
| 1194 | | 998 | 0.2217004114 | 1247 | | 1044 | 0.2203499828 |
| 1195 | | 999 | 0.2216742261 | 1248 | | 1045 | 0.2203252139 |
| 1196 | | 1000 | 0.2216480689 | 1249 | 204 | | 0.2204235192 |
| 1197 | | 1001 | 0.2216219397 | 1250 | | 1046 | 0.2202757438 |
| 1198 | | 1002 | 0.2215958385 | 1251 | | 1047 | 0.2202510510 |
| 1199 | | 1003 | 0.2215697652 | 1252 | | 1048 | 0.2202263835 |
| 1200 | | 1004 | 0.2215437197 | 1253 | | 1049 | 0.2202017411 |
| 1201 | 197 | | 0.2216478323 | 1254 | | 1050 | 0.2201771240 |
| 1202 | | 1005 | 0.2214917026 | 1255 | | 1051 | 0.2201525319 |
| 1203 | | 1006 | 0.2214657404 | 1256 | | 1052 | 0.2201279649 |
| 1204 | | 1007 | 0.2214398058 | 1257 | | 1053 | 0.2201034230 |
| 1205 | | 1008 | 0.2214138988 | 1258 | | 1054 | 0.2200789060 |
| 1206 | | 1009 | 0.2213880194 | 1259 | 205 | | 0.2201765851 |
| 1207 | | 1010 | 0.2213621674 | 1260 | | 1055 | 0.2200299382 |
| 1208 | | 1011 | 0.2213363429 | 1261 | | 1056 | 0.2200054959 |
| 1209 | | 1012 | 0.2213105457 | 1262 | | 1057 | 0.2199810783 |
| 1210 | | 1013 | 0.2212847759 | 1263 | | 1058 | 0.2199566855 |
| 1211 | | 1014 | 0.2212590333 | 1264 | | 1059 | 0.2199323174 |

## The Complete Table of Primes (first 1000 Primes, less than 7919)

| $m$ | $p_m$ | $t_m$ | $m$ | $p_m$ | $t_m$ |
|---|---|---|---|---|---|
| 1 | 2 | 3.399270106 | 46 | 199 | 0.2979842046 |
| 2 | 3 | 1.6963574120 | 47 | 211 | 0.2947079482 |
| 3 | 5 | 1.0924502700 | 48 | 223 | 0.2916766064 |
| 4 | 7 | 0.8754118490 | 49 | 227 | 0.2906839621 |
| 5 | 11 | 0.7012720514 | 50 | 229 | 0.2901750901 |
| 6 | 13 | 0.6459186318 | 51 | 233 | 0.2892205137 |
| 7 | 17 | 0.5804327797 | 52 | 239 | 0.2878510528 |
| 8 | 19 | 0.5540524507 | 53 | 241 | 0.2873799363 |
| 9 | 23 | 0.5179753118 | 54 | 251 | 0.2852491381 |
| 10 | 29 | 0.4810304852 | 55 | 257 | 0.2840114764 |
| 11 | 31 | 0.4697131298 | 56 | 263 | 0.2828128834 |
| 12 | 37 | 0.4458141449 | 57 | 269 | 0.2816512393 |
| 13 | 41 | 0.4324638276 | 58 | 271 | 0.2812519790 |
| 14 | 43 | 0.4259351409 | 59 | 277 | 0.2801374363 |
| 15 | 47 | 0.4154036691 | 60 | 281 | 0.2794035680 |
| 16 | 53 | 0.4023796406 | 61 | 283 | 0.2790291789 |
| 17 | 59 | 0.3913883733 | 62 | 293 | 0.2773117551 |
| 18 | 61 | 0.3876633904 | 63 | 307 | 0.2750451905 |
| 19 | 67 | 0.3787050658 | 64 | 311 | 0.2744063405 |
| 20 | 71 | 0.3732065227 | 65 | 313 | 0.2740801622 |
| 21 | 73 | 0.3704244073 | 66 | 317 | 0.2734382910 |
| 22 | 79 | 0.3635156856 | 67 | 331 | 0.2714163208 |
| 23 | 83 | 0.3592122065 | 68 | 337 | 0.2705646773 |
| 24 | 89 | 0.3534526181 | 69 | 347 | 0.2692025338 |
| 25 | 97 | 0.3466774192 | 70 | 349 | 0.2689215808 |
| 26 | 101 | 0.3434607769 | 71 | 353 | 0.2683849381 |
| 27 | 103 | 0.3418152806 | 72 | 359 | 0.2676041052 |
| 28 | 107 | 0.3388815946 | 73 | 367 | 0.2665954600 |
| 29 | 109 | 0.3373863987 | 74 | 373 | 0.2658541489 |
| 30 | 113 | 0.3346936830 | 75 | 379 | 0.2651287847 |
| 31 | 127 | 0.3266093677 | 76 | 383 | 0.2646488391 |
| 32 | 131 | 0.3244248382 | 77 | 389 | 0.2639489036 |
| 33 | 137 | 0.3213891550 | 78 | 397 | 0.2630425465 |
| 34 | 139 | 0.3203388907 | 79 | 401 | 0.2625917536 |
| 35 | 149 | 0.3158475580 | 80 | 409 | 0.2617212887 |
| 36 | 151 | 0.3149160280 | 81 | 419 | 0.2606677385 |
| 37 | 157 | 0.3124271024 | 82 | 421 | 0.2604512155 |
| 38 | 163 | 0.3100684523 | 83 | 431 | 0.2594371894 |
| 39 | 167 | 0.3085350152 | 84 | 433 | 0.2592289073 |
| 40 | 173 | 0.3063704960 | 85 | 439 | 0.2586346845 |
| 41 | 179 | 0.3043085611 | 86 | 443 | 0.2582409792 |
| 42 | 181 | 0.3035965569 | 87 | 449 | 0.2576646712 |
| 43 | 191 | 0.3004605526 | 88 | 457 | 0.2569154590 |
| 44 | 193 | 0.2998110688 | 89 | 461 | 0.2565424308 |
| 45 | 197 | 0.2986023209 | 90 | 463 | 0.2563528155 |

| $m$ | $p_m$ | $t_m$ | $m$ | $p_m$ | $t_m$ |
|---|---|---|---|---|---|
| 91 | 467 | 0.2559865423 | 139 | 797 | 0.2353362155 |
| 92 | 479 | 0.2549304335 | 140 | 809 | 0.2348095014 |
| 93 | 487 | 0.2542426555 | 141 | 811 | 0.2347201774 |
| 94 | 491 | 0.2539000843 | 142 | 821 | 0.2342898411 |
| 95 | 499 | 0.2532344206 | 143 | 823 | 0.2342022563 |
| 96 | 503 | 0.2529028203 | 144 | 827 | 0.2340308596 |
| 97 | 509 | 0.2524160123 | 145 | 829 | 0.2339441783 |
| 98 | 521 | 0.2514724805 | 146 | 839 | 0.2335259566 |
| 99 | 523 | 0.2513117462 | 147 | 853 | 0.2329522571 |
| 100 | 541 | 0.2499593894 | 148 | 857 | 0.2327886730 |
| 101 | 547 | 0.2495169078 | 149 | 859 | 0.2327059579 |
| 102 | 557 | 0.2487981832 | 150 | 863 | 0.2325439179 |
| 103 | 563 | 0.2483723239 | 151 | 877 | 0.2319906989 |
| 104 | 569 | 0.2479524562 | 152 | 881 | 0.2318329550 |
| 105 | 571 | 0.2478094900 | 153 | 883 | 0.2317532235 |
| 106 | 577 | 0.2473974860 | 154 | 887 | 0.2315969249 |
| 107 | 587 | 0.2467273388 | 155 | 907 | 0.2308383193 |
| 108 | 593 | 0.2463299992 | 156 | 911 | 0.2306873261 |
| 109 | 599 | 0.2459379565 | 157 | 919 | 0.2303902230 |
| 110 | 601 | 0.2458046101 | 158 | 929 | 0.2300240816 |
| 111 | 607 | 0.2454195260 | 159 | 937 | 0.2297343752 |
| 112 | 613 | 0.2450394510 | 160 | 941 | 0.2295896319 |
| 113 | 617 | 0.2447870228 | 161 | 947 | 0.2293751047 |
| 114 | 619 | 0.2446590445 | 162 | 953 | 0.2291623405 |
| 115 | 631 | 0.2439281459 | 163 | 967 | 0.2286753287 |
| 116 | 641 | 0.2433319846 | 164 | 971 | 0.2285343692 |
| 117 | 643 | 0.2432103042 | 165 | 977 | 0.2283303882 |
| 118 | 647 | 0.2429733659 | 166 | 983 | 0.2281260296 |
| 119 | 653 | 0.2426239452 | 167 | 991 | 0.2278567067 |
| 120 | 659 | 0.2422787374 | 168 | 997 | 0.2276560732 |
| 121 | 661 | 0.2421616210 | 169 | 1009 | 0.2272613448 |
| 122 | 673 | 0.2414905743 | 170 | 1013 | 0.2271298658 |
| 123 | 677 | 0.2412675044 | 171 | 1019 | 0.2269348412 |
| 124 | 683 | 0.2409382744 | 172 | 1021 | 0.2268689592 |
| 125 | 691 | 0.2405065322 | 173 | 1031 | 0.2265492299 |
| 126 | 701 | 0.2399769559 | 174 | 1033 | 0.2264843591 |
| 127 | 709 | 0.2395595552 | 175 | 1039 | 0.2262942275 |
| 128 | 719 | 0.2390473381 | 176 | 1049 | 0.2259816094 |
| 129 | 727 | 0.2386434483 | 177 | 1051 | 0.2259182039 |
| 130 | 733 | 0.2383433957 | 178 | 1061 | 0.2256101702 |
| 131 | 739 | 0.2380465552 | 179 | 1063 | 0.2255477087 |
| 132 | 743 | 0.2378492295 | 180 | 1069 | 0.2253644888 |
| 133 | 751 | 0.2374622673 | 181 | 1087 | 0.2248258330 |
| 134 | 757 | 0.2371719251 | 182 | 1091 | 0.2247064700 |
| 135 | 761 | 0.2369835318 | 183 | 1093 | 0.2246462370 |
| 136 | 769 | 0.2366084497 | 184 | 1097 | 0.2245277476 |
| 137 | 773 | 0.2364211939 | 185 | 1103 | 0.2243518169 |
| 138 | 787 | 0.2357835936 | 186 | 1109 | 0.2241771210 |

| $m$ | $p_m$ | $t_m$ | $m$ | $p_m$ | $t_m$ |
|---|---|---|---|---|---|
| 187 | 1117 | 0.2239465673 | 241 | 1523 | 0.2144583422 |
| 188 | 1123 | 0.2237746859 | 242 | 1531 | 0.2142802077 |
| 189 | 1129 | 0.2236039887 | 243 | 1543 | 0.2140519355 |
| 190 | 1151 | 0.2229916331 | 244 | 1549 | 0.2139382768 |
| 191 | 1153 | 0.2229354535 | 245 | 1553 | 0.2138625749 |
| 192 | 1163 | 0.2226619183 | 246 | 1559 | 0.2137498533 |
| 193 | 1171 | 0.2224449867 | 247 | 1567 | 0.2136006549 |
| 194 | 1181 | 0.2221768161 | 248 | 1571 | 0.2135260673 |
| 195 | 1187 | 0.2220167854 | 249 | 1579 | 0.2133783196 |
| 196 | 1193 | 0.2218577964 | 250 | 1583 | 0.2133044583 |
| 197 | 1201 | 0.2216478323 | 251 | 1597 | 0.2130495164 |
| 198 | 1213 | 0.2213368572 | 252 | 1601 | 0.2129767136 |
| 199 | 1217 | 0.2212332446 | 253 | 1607 | 0.2128682830 |
| 200 | 1223 | 0.2210792693 | 254 | 1609 | 0.2128318133 |
| 201 | 1229 | 0.2209262660 | 255 | 1613 | 0.2127597131 |
| 202 | 1231 | 0.2208746833 | 256 | 1619 | 0.2126523179 |
| 203 | 1237 | 0.2207229689 | 257 | 1621 | 0.2126162054 |
| 204 | 1249 | 0.2204235192 | 258 | 1627 | 0.2125094904 |
| 205 | 1259 | 0.2201765851 | 259 | 1637 | 0.2123331604 |
| 206 | 1277 | 0.2197392723 | 260 | 1657 | 0.2119851893 |
| 207 | 1279 | 0.2196901674 | 261 | 1663 | 0.2118814035 |
| 208 | 1283 | 0.2195933861 | 262 | 1667 | 0.2118122751 |
| 209 | 1289 | 0.2194494963 | 263 | 1669 | 0.2117774896 |
| 210 | 1291 | 0.2194010042 | 264 | 1693 | 0.2113707373 |
| 211 | 1297 | 0.2192582684 | 265 | 1697 | 0.2113031604 |
| 212 | 1301 | 0.2191632326 | 266 | 1699 | 0.2112691557 |
| 213 | 1303 | 0.2191153390 | 267 | 1709 | 0.2111022373 |
| 214 | 1307 | 0.2190208810 | 268 | 1721 | 0.2109036788 |
| 215 | 1319 | 0.2187417521 | 269 | 1723 | 0.2108702799 |
| 216 | 1321 | 0.2186947105 | 270 | 1733 | 0.2107063008 |
| 217 | 1327 | 0.2185561371 | 271 | 1741 | 0.2105758676 |
| 218 | 1361 | 0.2177902581 | 272 | 1747 | 0.2104784028 |
| 219 | 1367 | 0.2176568353 | 273 | 1753 | 0.2103813640 |
| 220 | 1373 | 0.2175241636 | 274 | 1759 | 0.2102847479 |
| 221 | 1381 | 0.2173487375 | 275 | 1777 | 0.2099984682 |
| 222 | 1399 | 0.2169599007 | 276 | 1783 | 0.2099034995 |
| 223 | 1409 | 0.2167462119 | 277 | 1787 | 0.2098402366 |
| 224 | 1423 | 0.2164506441 | 278 | 1789 | 0.2098084116 |
| 225 | 1427 | 0.2163662251 | 279 | 1801 | 0.2096210199 |
| 226 | 1429 | 0.2163236821 | 280 | 1811 | 0.2094659781 |
| 227 | 1433 | 0.2162397304 | 281 | 1823 | 0.2092814556 |
| 228 | 1439 | 0.2161148104 | 282 | 1831 | 0.2091591240 |
| 229 | 1447 | 0.2159495744 | 283 | 1847 | 0.2089169815 |
| 230 | 1451 | 0.2158669648 | 284 | 1861 | 0.2087072112 |
| 231 | 1453 | 0.2158253479 | 285 | 1867 | 0.2086176371 |
| 232 | 1459 | 0.2157026327 | 286 | 1871 | 0.2085579614 |
| 233 | 1471 | 0.2154599545 | 287 | 1873 | 0.2085279429 |
| 234 | 1481 | 0.2152595043 | 288 | 1877 | 0.2084685157 |
| 235 | 1483 | 0.2152189718 | 289 | 1879 | 0.2084386257 |
| 236 | 1487 | 0.2151389327 | 290 | 1889 | 0.2082916767 |
| 237 | 1489 | 0.2150986227 | 291 | 1901 | 0.2081167252 |
| 238 | 1493 | 0.2150190072 | 292 | 1907 | 0.2080295407 |
| 239 | 1499 | 0.2149004870 | 293 | 1913 | 0.2079427045 |
| 240 | 1511 | 0.2146660191 | 294 | 1931 | 0.2076851606 |

| $m$ | $p_m$ | $t_m$ | $m$ | $p_m$ | $t_m$ |
|---|---|---|---|---|---|
| 295 | 1933 | 0.2076563298 | 349 | 2351 | 0.2023033750 |
| 296 | 1949 | 0.2074301825 | 350 | 2357 | 0.2022371968 |
| 297 | 1951 | 0.2074016898 | 351 | 2371 | 0.2020832147 |
| 298 | 1973 | 0.2070951015 | 352 | 2377 | 0.2020177361 |
| 299 | 1979 | 0.2070119165 | 353 | 2381 | 0.2019742933 |
| 300 | 1987 | 0.2069016445 | 354 | 2383 | 0.2019530380 |
| 301 | 1993 | 0.2068192012 | 355 | 2389 | 0.2018876779 |
| 302 | 1997 | 0.2067642722 | 356 | 2393 | 0.2018445046 |
| 303 | 1999 | 0.2067366498 | 357 | 2399 | 0.2017797743 |
| 304 | 2003 | 0.2065549514 | 358 | 2411 | 0.2016506466 |
| 305 | 2011 | 0.2064469972 | 359 | 2417 | 0.2015865214 |
| 306 | 2017 | 0.2064921424 | 360 | 2423 | 0.2015225950 |
| 307 | 2027 | 0.2063577694 | 361 | 2437 | 0.2013738447 |
| 308 | 2029 | 0.2062063706 | 362 | 2441 | 0.2013317323 |
| 309 | 2039 | 0.2060735787 | 363 | 2447 | 0.2012685918 |
| 310 | 2053 | 0.2058888859 | 364 | 2459 | 0.2011426282 |
| 311 | 2063 | 0.2057580463 | 365 | 2467 | 0.2010591677 |
| 312 | 2069 | 0.2056800824 | 366 | 2473 | 0.2009971169 |
| 313 | 2081 | 0.2055246167 | 367 | 2477 | 0.2009555110 |
| 314 | 2083 | 0.2054991370 | 368 | 2503 | 0.2006873698 |
| 315 | 2087 | 0.2054478826 | 369 | 2521 | 0.2005038534 |
| 316 | 2089 | 0.2054224900 | 370 | 2531 | 0.2004024684 |
| 317 | 2099 | 0.2052944741 | 371 | 2539 | 0.2003222272 |
| 318 | 2111 | 0.2051417921 | 372 | 2543 | 0.2002822236 |
| 319 | 2113 | 0.2051167596 | 373 | 2549 | 0.2002222404 |
| 320 | 2129 | 0.2049149938 | 374 | 2551 | 0.2002024490 |
| 321 | 2131 | 0.2048902261 | 375 | 2557 | 0.2001426992 |
| 322 | 2137 | 0.2048153701 | 376 | 2579 | 0.1999244658 |
| 323 | 2141 | 0.2047657313 | 377 | 2591 | 0.1998065217 |
| 324 | 2143 | 0.2047411320 | 378 | 2593 | 0.1997871308 |
| 325 | 2153 | 0.2046171527 | 379 | 2609 | 0.1996310036 |
| 326 | 2161 | 0.2045185602 | 380 | 2617 | 0.1995535082 |
| 327 | 2179 | 0.2042979666 | 381 | 2621 | 0.1995149890 |
| 328 | 2203 | 0.2040072760 | 382 | 2633 | 0.1993994033 |
| 329 | 2207 | 0.2039595011 | 383 | 2647 | 0.1992653468 |
| 330 | 2213 | 0.2038878706 | 384 | 2657 | 0.1991702013 |
| 331 | 2221 | 0.2037926289 | 385 | 2659 | 0.1991514092 |
| 332 | 2237 | 0.2036031009 | 386 | 2663 | 0.1991136481 |
| 333 | 2239 | 0.2035798206 | 387 | 2671 | 0.1990381112 |
| 334 | 2243 | 0.2035330048 | 388 | 2677 | 0.1989817000 |
| 335 | 2251 | 0.2034393582 | 389 | 2683 | 0.1989254463 |
| 336 | 2267 | 0.2032529917 | 390 | 2687 | 0.1988879942 |
| 337 | 2269 | 0.2032300942 | 391 | 2689 | 0.1988695690 |
| 338 | 2273 | 0.2031840516 | 392 | 2693 | 0.1988323292 |
| 339 | 2281 | 0.2030919522 | 393 | 2699 | 0.1987764896 |
| 340 | 2287 | 0.2030232220 | 394 | 2707 | 0.1987022057 |
| 341 | 2293 | 0.2029547172 | 395 | 2711 | 0.1986652733 |
| 342 | 2297 | 0.2029092747 | 396 | 2713 | 0.1986469383 |
| 343 | 2309 | 0.2027729257 | 397 | 2719 | 0.1985916100 |
| 344 | 2311 | 0.2027505415 | 398 | 2729 | 0.1984995951 |
| 345 | 2333 | 0.2025028837 | 399 | 2731 | 0.1984814090 |
| 346 | 2339 | 0.2024360692 | 400 | 2741 | 0.1983898981 |
| 347 | 2341 | 0.2024140440 | 401 | 2749 | 0.1983170308 |
| 348 | 2347 | 0.2023475131 | 402 | 2753 | 0.1982807983 |

| $m$ | $p_m$ | $t_m$ | $m$ | $p_m$ | $t_m$ |
|---|---|---|---|---|---|
| 403 | 2767 | 0.1981454071 | 457 | 3229 | 0.1942776729 |
| 404 | 2777 | 0.1980375993 | 458 | 3251 | 0.1942627625 |
| 405 | 2789 | 0.1980195209 | 459 | 3253 | 0.1942331232 |
| 406 | 2791 | 0.1979657837 | 460 | 3257 | 0.1942182485 |
| 407 | 2797 | 0.1979299737 | 461 | 3259 | 0.1941299616 |
| 408 | 2801 | 0.1979119950 | 462 | 3271 | 0.1939257108 |
| 409 | 2803 | 0.1977701104 | 463 | 3299 | 0.1939110715 |
| 410 | 2819 | 0.1976467611 | 464 | 3301 | 0.1938675059 |
| 411 | 2833 | 0.1976115208 | 465 | 3307 | 0.1938240391 |
| 412 | 2837 | 0.1975588738 | 466 | 3313 | 0.1937806707 |
| 413 | 2843 | 0.1974889578 | 467 | 3319 | 0.1937517666 |
| 414 | 2851 | 0.1974366348 | 468 | 3323 | 0.1937085617 |
| 415 | 2857 | 0.1974017669 | 469 | 3329 | 0.1936940899 |
| 416 | 2861 | 0.1972462678 | 470 | 3331 | 0.1936081739 |
| 417 | 2879 | 0.1971774436 | 471 | 3343 | 0.1935795298 |
| 418 | 2887 | 0.1970917947 | 472 | 3347 | 0.1934941249 |
| 419 | 2897 | 0.1970405095 | 473 | 3359 | 0.1934798154 |
| 420 | 2903 | 0.1969893572 | 474 | 3361 | 0.1934089653 |
| 421 | 2909 | 0.1969214202 | 475 | 3371 | 0.1933947201 |
| 422 | 2917 | 0.1968026931 | 476 | 3373 | 0.1932820666 |
| 423 | 2927 | 0.1967359241 | 477 | 3389 | 0.1932679160 |
| 424 | 2939 | 0.1966188314 | 478 | 3391 | 0.1931560057 |
| 425 | 2953 | 0.1965853773 | 479 | 3407 | 0.1931141264 |
| 426 | 2957 | 0.1965353911 | 480 | 3413 | 0.1929754916 |
| 427 | 2963 | 0.1964855318 | 481 | 3433 | 0.1928652800 |
| 428 | 2969 | 0.1964688244 | 482 | 3449 | 0.1928103483 |
| 429 | 2971 | 0.1962386319 | 483 | 3457 | 0.1927828780 |
| 430 | 2999 | 0.1962221326 | 484 | 3461 | 0.1927690942 |
| 431 | 3001 | 0.1961405290 | 485 | 3463 | 0.1927416842 |
| 432 | 3011 | 0.1960754557 | 486 | 3467 | 0.1927279313 |
| 433 | 3019 | 0.1960429156 | 487 | 3469 | 0.1925785512 |
| 434 | 3023 | 0.1959298685 | 488 | 3491 | 0.1925244421 |
| 435 | 3037 | 0.1958975696 | 489 | 3499 | 0.1924436288 |
| 436 | 3041 | 0.1958332980 | 490 | 3511 | 0.1924032894 |
| 437 | 3049 | 0.1957373672 | 491 | 3517 | 0.1923363293 |
| 438 | 3061 | 0.1956894959 | 492 | 3527 | 0.1923228680 |
| 439 | 3067 | 0.1955942670 | 493 | 3529 | 0.1922960965 |
| 440 | 3079 | 0.1955625189 | 494 | 3533 | 0.1922560707 |
| 441 | 3083 | 0.1955150746 | 495 | 3539 | 0.1922426677 |
| 442 | 3089 | 0.1953581286 | 496 | 3541 | 0.1922027553 |
| 443 | 3109 | 0.1952800477 | 497 | 3547 | 0.1921365006 |
| 444 | 3119 | 0.1952643422 | 498 | 3557 | 0.1921231830 |
| 445 | 3121 | 0.1951402554 | 499 | 3559 | 0.1920440622 |
| 446 | 3137 | 0.1949403874 | 500 | 3571 | 0.1919783610 |
| 447 | 3163 | 0.1949096883 | 501 | 3581 | 0.1919651552 |
| 448 | 3167 | 0.1948942795 | 502 | 3583 | 0.1918997280 |
| 449 | 3169 | 0.1948028560 | 503 | 3593 | 0.1918085556 |
| 450 | 3181 | 0.1947572289 | 504 | 3607 | 0.1917695509 |
| 451 | 3187 | 0.1947268197 | 505 | 3613 | 0.1917435548 |
| 452 | 3191 | 0.1946361816 | 506 | 3617 | 0.1917046849 |
| 453 | 3203 | 0.1945909459 | 507 | 3623 | 0.1916530207 |
| 454 | 3209 | 0.1945308461 | 508 | 3631 | 0.1916143375 |
| 455 | 3217 | 0.1945007917 | 509 | 3637 | 0.1915757340 |
| 456 | 3221 | 0.1944409716 | 510 | 3643 | 0.1914733619 |

| $m$ | $p_m$ | $t_m$ | $m$ | $p_m$ | $t_m$ |
|---|---|---|---|---|---|
| **511** | **3659** | **0.1913969198** | **563** | **4091** | **0.1888892009** |
| **512** | **3671** | **0.1913841180** | **564** | **4093** | **0.1888558705** |
| **513** | **3673** | **0.1913586508** | **565** | **4099** | **0.1887894764** |
| **514** | **3677** | **0.1912700590** | **566** | **4111** | **0.1887013522** |
| **515** | **3691** | **0.1912321570** | **567** | **4127** | **0.1886902924** |
| **516** | **3697** | **0.1912068958** | **568** | **4129** | **0.1886682776** |
| **517** | **3701** | **0.1911565828** | **569** | **4133** | **0.1886353473** |
| **518** | **3709** | **0.1910939079** | **570** | **4139** | **0.1885588510** |
| **519** | **3719** | **0.1910438976** | **571** | **4153** | **0.1885369943** |
| **520** | **3727** | **0.1910064511** | **572** | **4157** | **0.1885260345** |
| **521** | **3733** | **0.1909690796** | **573** | **4159** | **0.1884283494** |
| **522** | **3739** | **0.1908329648** | **574** | **4177** | **0.1882989384** |
| **523** | **3761** | **0.1907959385** | **575** | **4201** | **0.1882452397** |
| **524** | **3767** | **0.1907835437** | **576** | **4211** | **0.1882130638** |
| **525** | **3769** | **0.1907221060** | **577** | **4217** | **0.1882022972** |
| **526** | **3779** | **0.1906364731** | **578** | **4219** | **0.1881488825** |
| **527** | **3793** | **0.1906120057** | **579** | **4229** | **0.1881381542** |
| **528** | **3797** | **0.1905754159** | **580** | **4231** | **0.1880849273** |
| **529** | **3803** | **0.1904662780** | **581** | **4241** | **0.1880742370** |
| **530** | **3821** | **0.1904541015** | **582** | **4243** | **0.1880211968** |
| **531** | **3823** | **0.1903937397** | **583** | **4253** | **0.1879894154** |
| **532** | **3833** | **0.1903096011** | **584** | **4259** | **0.1879787817** |
| **533** | **3847** | **0.1902855605** | **585** | **4261** | **0.1879260193** |
| **534** | **3851** | **0.1902735026** | **586** | **4271** | **0.1879154230** |
| **535** | **3853** | **0.1902137241** | **587** | **4273** | **0.1878628443** |
| **536** | **3863** | **0.1901303956** | **588** | **4283** | **0.1878313393** |
| **537** | **3877** | **0.1901065865** | **589** | **4289** | **0.1877894432** |
| **538** | **3881** | **0.1900591571** | **590** | **4297** | **0.1876333989** |
| **539** | **3889** | **0.1899530032** | **591** | **4327** | **0.1875816309** |
| **540** | **3907** | **0.1899294212** | **592** | **4337** | **0.1875712344** |
| **541** | **3911** | **0.1898941526** | **593** | **4339** | **0.1875196438** |
| **542** | **3917** | **0.1898823480** | **594** | **4349** | **0.1874784621** |
| **543** | **3919** | **0.1898588566** | **595** | **4357** | **0.1874476186** |
| **544** | **3923** | **0.1898237228** | **596** | **4363** | **0.1873963793** |
| **545** | **3929** | **0.1898119640** | **597** | **4373** | **0.1873045725** |
| **546** | **3931** | **0.1897420341** | **598** | **4391** | **0.1872740244** |
| **547** | **3943** | **0.1897187210** | **599** | **4397** | **0.1872131568** |
| **548** | **3947** | **0.1896029636** | **600** | **4409** | **0.1871524939** |
| **549** | **3967** | **0.1894764731** | **601** | **4421** | **0.1871423425** |
| **550** | **3989** | **0.1894078010** | **602** | **4423** | **0.1870518167** |
| **551** | **4001** | **0.1893963045** | **603** | **4441** | **0.1870216941** |
| **552** | **4003** | **0.1893734244** | **604** | **4447** | **0.1870016164** |
| **553** | **4007** | **0.1893392029** | **605** | **4451** | **0.1869715778** |
| **554** | **4013** | **0.1893050451** | **606** | **4457** | **0.1869415895** |
| **555** | **4019** | **0.1892936137** | **607** | **4463** | **0.1868520655** |
| **556** | **4021** | **0.1892595408** | **608** | **4481** | **0.1868420829** |
| **557** | **4027** | **0.1891353795** | **609** | **4483** | **0.1867925356** |
| **558** | **4049** | **0.1891240534** | **610** | **4493** | **0.1867234263** |
| **559** | **4051** | **0.1890902933** | **611** | **4507** | **0.1866938493** |
| **560** | **4057** | **0.1890007142** | **612** | **4513** | **0.1866741352** |
| **561** | **4073** | **0.1889671801** | **613** | **4517** | **0.1866642518** |
| **562** | **4079** | **0.1889003821** | **614** | **4519** | **0.1866445706** |

| $m$ | $p_m$ | $t_m$ | $m$ | $p_m$ | $t_m$ |
|---|---|---|---|---|---|
| **615** | **4523** | **0.1865272779** | **668** | **4993** | **0.1844494303** |
| **616** | **4547** | **0.1865174753** | **669** | **4999** | **0.1844320621** |
| **617** | **4549** | **0.1864591135** | **670** | **5003** | **0.1844060696** |
| **618** | **4561** | **0.1864299698** | **671** | **5009** | **0.1843973774** |
| **619** | **4567** | **0.1863525974** | **672** | **5011** | **0.1843542003** |
| **620** | **4583** | **0.1863140027** | **673** | **5021** | **0.1843455340** |
| **621** | **4591** | **0.1862850946** | **674** | **5023** | **0.1842767376** |
| **622** | **4597** | **0.1862562332** | **675** | **5039** | **0.1842253035** |
| **623** | **4603** | **0.1861700609** | **676** | **5051** | **0.1841910799** |
| **624** | **4621** | **0.1860938037** | **677** | **5059** | **0.1841143877** |
| **625** | **4637** | **0.1860842366** | **678** | **5077** | **0.1840973489** |
| **626** | **4639** | **0.1860651840** | **679** | **5081** | **0.1840718486** |
| **627** | **4643** | **0.1860366762** | **680** | **5087** | **0.1840210124** |
| **628** | **4649** | **0.1860271402** | **681** | **5099** | **0.1840125100** |
| **629** | **4651** | **0.1859986934** | **682** | **5101** | **0.1839871332** |
| **630** | **4657** | **0.1859702919** | **683** | **5107** | **0.1839617934** |
| **631** | **4663** | **0.1859230996** | **684** | **5113** | **0.1839364903** |
| **632** | **4673** | **0.1858948185** | **685** | **5119** | **0.1838190825** |
| **633** | **4679** | **0.1858384547** | **686** | **5147** | **0.1837939854** |
| **634** | **4691** | **0.1857822689** | **687** | **5153** | **0.1837356345** |
| **635** | **4703** | **0.1856983534** | **688** | **5167** | **0.1837189618** |
| **636** | **4721** | **0.1856889978** | **689** | **5171** | **0.1836857157** |
| **637** | **4723** | **0.1856610865** | **690** | **5179** | **0.1836442599** |
| **638** | **4729** | **0.1856424826** | **691** | **5189** | **0.1836111561** |
| **639** | **4733** | **0.1855592233** | **692** | **5197** | **0.1835616435** |
| **640** | **4751** | **0.1855223097** | **693** | **5209** | **0.1834876632** |
| **641** | **4759** | **0.1854121512** | **694** | **5227** | **0.1834712266** |
| **642** | **4783** | **0.1853938064** | **695** | **5231** | **0.1834629889** |
| **643** | **4787** | **0.1853846105** | **696** | **5233** | **0.1834465758** |
| **644** | **4789** | **0.1853662946** | **697** | **5237** | **0.1833486704** |
| **645** | **4793** | **0.1853388867** | **698** | **5261** | **0.1832998991** |
| **646** | **4799** | **0.1853297199** | **699** | **5273** | **0.1832755399** |
| **647** | **4801** | **0.1852751201** | **700** | **5279** | **0.1832673949** |
| **648** | **4813** | **0.1852569175** | **701** | **5281** | **0.1832027166** |
| **649** | **4817** | **0.1851935053** | **702** | **5297** | **0.1831784937** |
| **650** | **4831** | **0.1850584514** | **703** | **5303** | **0.1831543046** |
| **651** | **4861** | **0.1850136223** | **704** | **5309** | **0.1830980594** |
| **652** | **4871** | **0.1849867558** | **705** | **5323** | **0.1830579820** |
| **653** | **4877** | **0.1849332042** | **706** | **5333** | **0.1830020483** |
| **654** | **4889** | **0.1848709420** | **707** | **5347** | **0.1829860660** |
| **655** | **4903** | **0.1848442920** | **708** | **5351** | **0.1828669784** |
| **656** | **4909** | **0.1848000036** | **709** | **5381** | **0.1828432207** |
| **657** | **4919** | **0.1847470157** | **710** | **5387** | **0.1828194958** |
| **658** | **4931** | **0.1847381516** | **711** | **5393** | **0.1827958034** |
| **659** | **4933** | **0.1847204948** | **712** | **5399** | **0.1827642800** |
| **660** | **4937** | **0.1846940717** | **713** | **5407** | **0.1827406633** |
| **661** | **4943** | **0.1846589217** | **714** | **5413** | **0.1827249212** |
| **662** | **4951** | **0.1846325910** | **715** | **5417** | **0.1827170322** |
| **663** | **4957** | **0.1845888318** | **716** | **5419** | **0.1826700072** |
| **664** | **4967** | **0.1845800478** | **717** | **5431** | **0.1826465195** |
| **665** | **4969** | **0.1845625496** | **718** | **5437** | **0.1826308634** |
| **666** | **4973** | **0.1845015831** | **719** | **5441** | **0.1826230176** |
| **667** | **4987** | **0.1844754873** | **720** | **5443** | **0.1825995940** |

| $m$ | $p_m$ | $t_m$ | $m$ | $p_m$ | $t_m$ |
|---|---|---|---|---|---|
| **721** | **5449** | **0.1825141009** | **772** | **5867** | **0.1810355112** |
| **722** | **5471** | **0.1824908248** | **773** | **5869** | **0.1809999754** |
| **723** | **5477** | **0.1824830428** | **774** | **5879** | **0.1809928466** |
| **724** | **5479** | **0.1824675350** | **775** | **5881** | **0.1809361903** |
| **725** | **5483** | **0.1823980799** | **776** | **5897** | **0.1809149693** |
| **726** | **5501** | **0.1823903399** | **777** | **5903** | **0.1808445118** |
| **727** | **5503** | **0.1823749155** | **778** | **5923** | **0.1808304253** |
| **728** | **5507** | **0.1823288146** | **779** | **5927** | **0.1807883111** |
| **729** | **5519** | **0.1823211061** | **780** | **5939** | **0.1807393163** |
| **730** | **5521** | **0.1822980901** | **781** | **5953** | **0.1806417875** |
| **731** | **5527** | **0.1822827485** | **782** | **5981** | **0.1806209323** |
| **732** | **5531** | **0.1821836026** | **783** | **5987** | **0.1805516867** |
| **733** | **5557** | **0.1821607701** | **784** | **6007** | **0.1805378421** |
| **734** | **5563** | **0.1821379681** | **785** | **6011** | **0.1804758101** |
| **735** | **5569** | **0.1821227690** | **786** | **6029** | **0.1804482930** |
| **736** | **5573** | **0.1820924549** | **787** | **6037** | **0.1804276754** |
| **737** | **5581** | **0.1820546482** | **788** | **6043** | **0.1804139323** |
| **738** | **5591** | **0.1819343191** | **789** | **6047** | **0.1803933568** |
| **739** | **5623** | **0.1818744487** | **790** | **6053** | **0.1803454931** |
| **740** | **5639** | **0.1818669421** | **791** | **6067** | **0.1803250010** |
| **741** | **5641** | **0.1818445277** | **792** | **6073** | **0.1803045339** |
| **742** | **5647** | **0.1818295870** | **793** | **6079** | **0.1802705012** |
| **743** | **5651** | **0.1818221004** | **794** | **6089** | **0.1802636742** |
| **744** | **5653** | **0.1818071794** | **795** | **6091** | **0.1802297240** |
| **745** | **5657** | **0.1817997028** | **796** | **6101** | **0.1801890810** |
| **746** | **5659** | **0.1817625367** | **797** | **6113** | **0.1801620283** |
| **747** | **5669** | **0.1817106565** | **798** | **6121** | **0.1801282824** |
| **748** | **5683** | **0.1816884464** | **799** | **6131** | **0.1801215132** |
| **749** | **5689** | **0.1816736419** | **800** | **6133** | **0.1800878485** |
| **750** | **5693** | **0.1816441129** | **801** | **6143** | **0.1800609582** |
| **751** | **5701** | **0.1816072838** | **802** | **6151** | **0.1800207208** |
| **752** | **5711** | **0.1815852079** | **803** | **6163** | **0.1799872573** |
| **753** | **5717** | **0.1815119229** | **804** | **6173** | **0.1799072646** |
| **754** | **5737** | **0.1814972709** | **805** | **6197** | **0.1799005841** |
| **755** | **5741** | **0.1814899294** | **806** | **6199** | **0.1798872655** |
| **756** | **5743** | **0.1814680050** | **807** | **6203** | **0.1798606946** |
| **757** | **5749** | **0.1813589657** | **808** | **6211** | **0.1798407855** |
| **758** | **5779** | **0.1813444448** | **809** | **6217** | **0.1798275145** |
| **759** | **5783** | **0.1813154805** | **810** | **6221** | **0.1798010381** |
| **760** | **5791** | **0.1812793543** | **811** | **6229** | **0.1797416611** |
| **761** | **5801** | **0.1812576994** | **812** | **6247** | **0.1797087497** |
| **762** | **5807** | **0.1812360721** | **813** | **6257** | **0.1796890203** |
| **763** | **5813** | **0.1812072917** | **814** | **6263** | **0.1796693143** |
| **764** | **5821** | **0.1811857286** | **815** | **6269** | **0.1796627285** |
| **765** | **5827** | **0.1811427235** | **816** | **6271** | **0.1796430534** |
| **766** | **5839** | **0.1811283865** | **817** | **6277** | **0.1796103353** |
| **767** | **5843** | **0.1811069234** | **818** | **6287** | **0.1795711644** |
| **768** | **5849** | **0.1810997491** | **819** | **6299** | **0.1795646173** |
| **769** | **5851** | **0.1810783222** | **820** | **6301** | **0.1795320522** |
| **770** | **5857** | **0.1810640398** | **821** | **6311** | **0.1795125304** |
| **771** | **5861** | **0.1810426580** | **822** | **6317** | **0.1794930315** |

| $m$ | $p_m$ | $t_m$ | $m$ | $p_m$ | $t_m$ |
|---|---|---|---|---|---|
| **823** | **6323** | **0.1794735553** | **874** | **6791** | **0.1780334063** |
| **824** | **6329** | **0.1794476333** | **875** | **6793** | **0.1780037117** |
| **825** | **6337** | **0.1794282100** | **876** | **6803** | **0.1779445108** |
| **826** | **6343** | **0.1793959097** | **877** | **6823** | **0.1779326739** |
| **827** | **6353** | **0.1793765466** | **878** | **6827** | **0.1779267449** |
| **828** | **6359** | **0.1793700758** | **879** | **6829** | **0.1779149209** |
| **829** | **6361** | **0.1793507428** | **880** | **6833** | **0.1778913260** |
| **830** | **6367** | **0.1793314321** | **881** | **6841** | **0.1778442650** |
| **831** | **6373** | **0.1793121439** | **882** | **6857** | **0.1778266348** |
| **832** | **6379** | **0.1792800674** | **883** | **6863** | **0.1778090235** |
| **833** | **6389** | **0.1792544443** | **884** | **6869** | **0.1778031392** |
| **834** | **6397** | **0.1791778735** | **885** | **6871** | **0.1777680129** |
| **835** | **6421** | **0.1791587622** | **886** | **6883** | **0.1777213035** |
| **836** | **6427** | **0.1790889577** | **887** | **6899** | **0.1776979852** |
| **837** | **6449** | **0.1790825980** | **888** | **6907** | **0.1776863250** |
| **838** | **6451** | **0.1790257184** | **889** | **6911** | **0.1776688637** |
| **839** | **6469** | **0.1790130807** | **890** | **6917** | **0.1775819418** |
| **840** | **6473** | **0.1789878652** | **891** | **6947** | **0.1775761385** |
| **841** | **6481** | **0.1789564072** | **892** | **6949** | **0.1775472592** |
| **842** | **6491** | **0.1788624507** | **893** | **6959** | **0.1775414683** |
| **843** | **6521** | **0.1788374626** | **894** | **6961** | **0.1775241607** |
| **844** | **6529** | **0.1787814151** | **895** | **6967** | **0.1775126237** |
| **845** | **6547** | **0.1787689620** | **896** | **6971** | **0.1774953467** |
| **846** | **6551** | **0.1787627238** | **897** | **6977** | **0.1774780879** |
| **847** | **6553** | **0.1787316894** | **898** | **6983** | **0.1774551131** |
| **848** | **6563** | **0.1787130845** | **899** | **6991** | **0.1774378967** |
| **849** | **6569** | **0.1787068674** | **900** | **6997** | **0.1774264206** |
| **850** | **6571** | **0.1786882905** | **901** | **7001** | **0.1773920921** |
| **851** | **6577** | **0.1786759075** | **902** | **7013** | **0.1773749420** |
| **852** | **6581** | **0.1786204026** | **903** | **7019** | **0.1773521118** |
| **853** | **6599** | **0.1785957768** | **904** | **7027** | **0.1773179389** |
| **854** | **6607** | **0.1785589215** | **905** | **7039** | **0.1773065467** |
| **855** | **6619** | **0.1785038070** | **906** | **7043** | **0.1772668005** |
| **856** | **6637** | **0.1784549667** | **907** | **7057** | **0.1772328055** |
| **857** | **6653** | **0.1784366703** | **908** | **7069** | **0.1772045261** |
| **858** | **6659** | **0.1784305565** | **909** | **7079** | **0.1771368901** |
| **859** | **6661** | **0.1783940674** | **910** | **7103** | **0.1771200060** |
| **860** | **6673** | **0.1783758384** | **911** | **7109** | **0.1770863152** |
| **861** | **6679** | **0.1783455208** | **912** | **7121** | **0.1770694834** |
| **862** | **6689** | **0.1783394408** | **913** | **7127** | **0.1770638600** |
| **863** | **6691** | **0.1783091900** | **914** | **7129** | **0.1770023794** |
| **864** | **6701** | **0.1783031235** | **915** | **7151** | **0.1769800644** |
| **865** | **6703** | **0.1782849948** | **916** | **7159** | **0.1769299982** |
| **866** | **6709** | **0.1782548436** | **917** | **7177** | **0.1769022391** |
| **867** | **6719** | **0.1782127357** | **918** | **7187** | **0.1768855965** |
| **868** | **6733** | **0.1782007043** | **919** | **7193** | **0.1768468627** |
| **869** | **6737** | **0.1781288414** | **920** | **7207** | **0.1768357952** |
| **870** | **6761** | **0.1781228410** | **921** | **7211** | **0.1768302521** |
| **871** | **6763** | **0.1780751139** | **922** | **7213** | **0.1768136831** |
| **872** | **6779** | **0.1780691331** | **923** | **7219** | **0.1767861218** |
| **873** | **6781** | **0.1780393740** | **924** | **7229** | **0.1767641015** |

| $m$ | $p_m$ | $t_m$ | $m$ | $p_m$ | $t_m$ |
|---|---|---|---|---|---|
| **925** | **7237** | **0.1767475999** | **963** | **7583** | **0.1758238357** |
| **926** | **7243** | **0.1767365970** | **964** | **7589** | **0.1758186299** |
| **927** | **7247** | **0.1767201263** | **965** | **7591** | **0.1757875389** |
| **928** | **7253** | **0.1766381034** | **966** | **7603** | **0.1757771742** |
| **929** | **7283** | **0.1765999547** | **967** | **7607** | **0.1757410041** |
| **930** | **7297** | **0.1765727537** | **968** | **7621** | **0.1756946252** |
| **931** | **7307** | **0.1765672999** | **969** | **7639** | **0.1756843202** |
| **932** | **7309** | **0.1765347345** | **970** | **7643** | **0.1756688859** |
| **933** | **7321** | **0.1765076427** | **971** | **7649** | **0.1756175948** |
| **934** | **7331** | **0.1765022108** | **972** | **7669** | **0.1756073392** |
| **935** | **7333** | **0.1764589865** | **973** | **7673** | **0.1755868688** |
| **936** | **7349** | **0.1764535711** | **974** | **7681** | **0.1755715279** |
| **937** | **7351** | **0.1764051004** | **975** | **7687** | **0.1755613018** |
| **938** | **7369** | **0.1763407056** | **976** | **7691** | **0.1755408901** |
| **939** | **7393** | **0.1762925714** | **977** | **7699** | **0.1755306835** |
| **940** | **7411** | **0.1762765429** | **978** | **7703** | **0.1754950639** |
| **941** | **7417** | **0.1762339166** | **979** | **7717** | **0.1754798107** |
| **942** | **7433** | **0.1761860990** | **980** | **7723** | **0.1754696429** |
| **943** | **7451** | **0.1761701758** | **981** | **7727** | **0.1754341586** |
| **944** | **7457** | **0.1761648563** | **982** | **7741** | **0.1754038029** |
| **945** | **7459** | **0.1761172424** | **983** | **7753** | **0.1753936833** |
| **946** | **7477** | **0.1761066628** | **984** | **7757** | **0.1753886156** |
| **947** | **7481** | **0.1760908179** | **985** | **7759** | **0.1753130797** |
| **948** | **7487** | **0.1760855245** | **986** | **7789** | **0.1753030173** |
| **949** | **7489** | **0.1760591747** | **987** | **7793** | **0.1752428778** |
| **950** | **7499** | **0.1760381212** | **988** | **7817** | **0.1752278629** |
| **951** | **7507** | **0.1760118487** | **989** | **7823** | **0.1752128621** |
| **952** | **7517** | **0.1759960967** | **990** | **7829** | **0.1751829230** |
| **953** | **7523** | **0.1759803601** | **991** | **7841** | **0.1751530399** |
| **954** | **7529** | **0.1759594094** | **992** | **7853** | **0.1751182502** |
| **955** | **7537** | **0.1759489331** | **993** | **7867** | **0.1751033520** |
| **956** | **7541** | **0.1759332426** | **994** | **7873** | **0.1750934208** |
| **957** | **7547** | **0.1759280010** | **995** | **7877** | **0.1750884475** |
| **958** | **7549** | **0.1759019074** | **996** | **7879** | **0.1750785257** |
| **959** | **7559** | **0.1758966761** | **997** | **7883** | **0.1750340236** |
| **960** | **7561** | **0.1758654340** | **998** | **7901** | **0.1750192038** |
| **961** | **7573** | **0.1758550188** | **999** | **7907** | **0.1749896254** |
| **962** | **7577** | **0.1758394197** | **1000** | **7919** | **0.1749699303** |

| $m$ | $p_m$ | $t_m$ | $t_0$ |
|---|---|---|---|
| 1 | 2 | 3.399270106 | 2.2661880071 |
| 2 | 3 | 1.6963574120 | 1.4298004340 |
| 3 | 5 | 1.0924502700 | 0.9759906329 |
| 4 | 7 | 0.8754118490 | 0.8072296285 |
| 5 | 11 | 0.7012720514 | 0.6550729486 |
| 6 | 13 | 0.6459186318 | 0.6124083200 |
| 7 | 17 | 0.5804327797 | 0.5544221829 |
| 8 | 19 | 0.5540524507 | 0.5334789880 |
| 9 | 23 | 0.5179753118 | 0.5009724843 |
| 10 | 29 | 0.4810304852 | 0.4664859894 |
| 11 | 31 | 0.4697131298 | 0.4574263773 |
| 12 | 37 | 0.4458141449 | 0.4350130257 |
| 13 | 41 | 0.4324638276 | 0.4229879745 |
| 14 | 43 | 0.4259351409 | 0.4176316810 |
| 15 | 47 | 0.4154036691 | 0.4079834046 |
| 16 | 53 | 0.4023796406 | 0.3956374899 |
| 17 | 59 | 0.3913883733 | 0.3852316131 |
| 18 | 61 | 0.3876633904 | 0.3821076439 |
| 19 | 67 | 0.3787050658 | 0.3735817262 |
| 20 | 71 | 0.3732065227 | 0.3684997166 |
| 21 | 73 | 0.3704244073 | 0.3661137807 |
| 22 | 79 | 0.3635156856 | 0.3594953824 |
| 23 | 83 | 0.3592122065 | 0.3554770281 |
| 24 | 89 | 0.3534526181 | 0.3499495609 |
| 25 | 97 | 0.3466774192 | 0.3433651512 |
| 26 | 101 | 0.3434607769 | 0.3403586799 |
| 27 | 103 | 0.3418152806 | 0.3389187007 |
| 28 | 107 | 0.3388815946 | 0.3361553317 |
| 29 | 109 | 0.3373863987 | 0.3348283620 |
| 30 | 113 | 0.3346936830 | 0.3322757487 |
| 31 | 127 | 0.3266093677 | 0.3242641745 |
| 32 | 131 | 0.3244248382 | 0.3222015895 |
| 33 | 137 | 0.3213891550 | 0.3192687839 |
| 34 | 139 | 0.3203388907 | 0.3183310618 |
| 35 | 149 | 0.3158475580 | 0.3139115072 |
| 36 | 151 | 0.3149160280 | 0.3130772805 |
| 37 | 157 | 0.3124271024 | 0.3106645498 |
| 38 | 163 | 0.3100684523 | 0.3083771808 |
| 39 | 167 | 0.3085350152 | 0.3069164177 |
| 40 | 173 | 0.3063704960 | 0.3048141752 |
| 41 | 179 | 0.3043085611 | 0.3028107771 |
| 42 | 181 | 0.3035965569 | 0.3021635518 |
| 43 | 191 | 0.3004605526 | 0.2990697930 |
| 44 | 193 | 0.2998110688 | 0.2984778261 |
| 45 | 197 | 0.2986023209 | 0.2973189011 |
| 46 | 199 | 0.2979842046 | 0.2967515340 |
| 47 | 211 | 0.2947079482 | 0.2935048515 |
| 48 | 223 | 0.2916766064 | 0.2905023907 |
| 49 | 227 | 0.2906839621 | 0.2895503777 |
| 50 | 229 | 0.2901750901 | 0.2890829391 |
| 51 | 233 | 0.2892205137 | 0.2881646021 |
| 52 | 239 | 0.2878510528 | 0.2868267655 |

| $m$ | $p_m$ | $t_m$ | $t_0$ |
|---|---|---|---|
| 53 | 241 | 0.2873799363 | 0.2863909723 |
| 54 | 251 | 0.2852491381 | 0.2842837219 |
| 55 | 257 | 0.2840114764 | 0.2830734888 |
| 56 | 263 | 0.2828128834 | 0.2819010957 |
| 57 | 269 | 0.2816512393 | 0.2807644971 |
| 58 | 271 | 0.2812519790 | 0.2803932542 |
| 59 | 277 | 0.2801374363 | 0.2793014647 |
| 60 | 281 | 0.2794035680 | 0.2785912591 |
| 61 | 283 | 0.2790291789 | 0.2782412722 |
| 62 | 293 | 0.2773117551 | 0.2765402453 |
| 63 | 307 | 0.2750451905 | 0.2742863781 |
| 64 | 311 | 0.2744063405 | 0.2736677700 |
| 65 | 313 | 0.2740801622 | 0.2733624751 |
| 66 | 317 | 0.2734382910 | 0.2727597012 |
| 67 | 331 | 0.2714163208 | 0.2707280730 |
| 68 | 337 | 0.2705646773 | 0.2698924303 |
| 69 | 347 | 0.2692025338 | 0.2685431885 |
| 70 | 349 | 0.2689215808 | 0.2682795955 |
| 71 | 353 | 0.2683849381 | 0.2677584386 |
| 72 | 359 | 0.2676041052 | 0.2669913738 |
| 73 | 367 | 0.2665954600 | 0.2659949326 |
| 74 | 373 | 0.2658541489 | 0.2652664907 |
| 75 | 379 | 0.2651287847 | 0.2645535577 |
| 76 | 383 | 0.2646488391 | 0.2640865983 |
| 77 | 389 | 0.2639489036 | 0.2633982440 |
| 78 | 397 | 0.2630425465 | 0.2625021814 |
| 79 | 401 | 0.2625917536 | 0.2620631357 |
| 80 | 409 | 0.2617212887 | 0.2612023160 |
| 81 | 419 | 0.2606677385 | 0.2601573216 |
| 82 | 421 | 0.2604512155 | 0.2599523039 |
| 83 | 431 | 0.2594371894 | 0.2589463134 |
| 84 | 433 | 0.2592289073 | 0.2587488370 |
| 85 | 439 | 0.2586346845 | 0.2581636089 |
| 86 | 443 | 0.2582409792 | 0.2577793290 |
| 87 | 449 | 0.2576646712 | 0.2572114688 |
| 88 | 457 | 0.2569154590 | 0.2564698004 |
| 89 | 461 | 0.2565424308 | 0.2561053947 |
| 90 | 463 | 0.2563528155 | 0.2559247606 |
| 91 | 467 | 0.2559865423 | 0.2555665766 |
| 92 | 479 | 0.2549304335 | 0.2545159640 |
| 93 | 487 | 0.2542426555 | 0.2538347258 |
| 94 | 491 | 0.2539000843 | 0.2534996353 |
| 95 | 499 | 0.2532344206 | 0.2528401612 |
| 96 | 503 | 0.2529028203 | 0.2525156430 |
| 97 | 509 | 0.2524160123 | 0.2520352071 |
| 98 | 521 | 0.2514724805 | 0.2510964019 |
| 99 | 523 | 0.2513117462 | 0.2509427085 |
| 100 | 541 | 0.2499593894 | 0.2495934654 |
| 101 | 547 | 0.2495169078 | 0.2491568059 |
| 102 | 557 | 0.2487981832 | 0.2484428815 |
| 103 | 563 | 0.2483723239 | 0.2480225758 |
| 104 | 569 | 0.2479524562 | 0.2476081230 |
| 105 | 571 | 0.2478094900 | 0.2474712476 |
| 106 | 577 | 0.2473974860 | 0.2470643751 |
| 107 | 587 | 0.2467273388 | 0.2463984645 |
| 108 | 593 | 0.2463299992 | 0.2460060303 |
| 109 | 599 | 0.2459379565 | 0.2456187761 |
| 110 | 601 | 0.2458046101 | 0.2454908215 |
| 111 | 607 | 0.2454195260 | 0.2451102863 |
| 112 | 613 | 0.2450394510 | 0.2447346548 |
| 113 | 617 | 0.2447870228 | 0.2444869025 |
| 114 | 619 | 0.2446590445 | 0.2443638149 |
| 115 | 631 | 0.2439281459 | 0.2436360800 |
| 116 | 641 | 0.2433319846 | 0.2430433494 |
| 117 | 643 | 0.2432103042 | 0.2429262557 |
| 118 | 647 | 0.2429733659 | 0.2426934925 |
| 119 | 653 | 0.2426239452 | 0.2423478569 |
| 120 | 659 | 0.2422787374 | 0.2420063524 |
| 121 | 661 | 0.2421616210 | 0.2418934203 |
| 122 | 673 | 0.2414905743 | 0.2412250859 |
| 123 | 677 | 0.2412675044 | 0.2410057613 |
| 124 | 683 | 0.2409382744 | 0.2406799301 |
| 125 | 691 | 0.2405065322 | 0.2402512594 |
| 126 | 701 | 0.2399769559 | 0.2397244466 |
| 127 | 709 | 0.2395595552 | 0.2393100093 |
| 128 | 719 | 0.2390473381 | 0.2388004613 |
| 129 | 727 | 0.2386434483 | 0.2383994316 |
| 130 | 733 | 0.2383433957 | 0.2381024155 |
| 131 | 739 | 0.2380465552 | 0.2378085509 |
| 132 | 743 | 0.2378492295 | 0.2376143626 |
| 133 | 751 | 0.2374622673 | 0.2372300399 |
| 134 | 757 | 0.2371719251 | 0.2369452792 |
| 135 | 761 | 0.2369835318 | 0.2367570659 |
| 136 | 769 | 0.2366084497 | 0.2363844730 |
| 137 | 773 | 0.2364211939 | 0.2362000625 |
| 138 | 787 | 0.2357835936 | 0.2355642711 |
| 139 | 797 | 0.2353362155 | 0.2351190682 |
| 140 | 809 | 0.2348095014 | 0.2345943100 |
| 141 | 811 | 0.2347201774 | 0.2345078331 |
| 142 | 821 | 0.2342898411 | 0.2340795645 |
| 143 | 823 | 0.2342022563 | 0.2339947230 |
| 144 | 827 | 0.2340308596 | 0.2338258402 |
| 145 | 829 | 0.2339441783 | 0.2337417957 |
| 146 | 839 | 0.2335259566 | 0.2333254848 |
| 147 | 853 | 0.2329522571 | 0.2327533402 |
| 148 | 857 | 0.2327886730 | 0.2325921028 |
| 149 | 859 | 0.2327059579 | 0.2325118495 |
| 150 | 863 | 0.2325439179 | 0.2323520669 |
| 151 | 877 | 0.2319906989 | 0.2318002960 |
| 152 | 881 | 0.2318329550 | 0.2316447394 |
| 153 | 883 | 0.2317532235 | 0.2315673037 |
| 154 | 887 | 0.2315969249 | 0.2314131114 |
| 155 | 907 | 0.2308383193 | 0.2306554288 |
| 156 | 911 | 0.2306873261 | 0.2305064845 |
| 157 | 919 | 0.2303902230 | 0.2302111180 |
| 158 | 929 | 0.2300240816 | 0.2298465516 |
| 159 | 937 | 0.2297343752 | 0.2295585316 |
| 160 | 941 | 0.2295896319 | 0.2294157107 |
| 161 | 947 | 0.2293751047 | 0.2292029438 |
| 162 | 953 | 0.2291623405 | 0.2289919112 |
| 163 | 967 | 0.2286753287 | 0.2285061055 |
| 164 | 971 | 0.2285343692 | 0.2283689694 |
| 165 | 977 | 0.2283303882 | 0.2281646267 |
| 166 | 983 | 0.2281260296 | 0.2279618970 |
| 167 | 991 | 0.2278567067 | 0.2276940608 |
| 168 | 997 | 0.2276560732 | 0.2274950073 |
| 169 | 1009 | 0.2272613448 | 0.2271014958 |
| 170 | 1013 | 0.2271298658 | 0.2269716639 |
| 171 | 1019 | 0.2269348412 | 0.2267781502 |
| 172 | 1021 | 0.2268689592 | 0.2267139717 |
| 173 | 1031 | 0.2265492299 | 0.2263954916 |
| 174 | 1033 | 0.2264843591 | 0.2263322730 |
| 175 | 1039 | 0.2262942275 | 0.2261435593 |
| 176 | 1049 | 0.2259816094 | 0.2258321347 |
| 177 | 1051 | 0.2259182039 | 0.2257703084 |
| 178 | 1061 | 0.2256101702 | 0.2254634327 |
| 179 | 1063 | 0.2255477087 | 0.2254025042 |
| 180 | 1069 | 0.2253644888 | 0.2252206003 |
| 181 | 1087 | 0.2248258330 | 0.2246826753 |
| 182 | 1091 | 0.2247064700 | 0.2245646912 |
| 183 | 1093 | 0.2246462370 | 0.2245059075 |
| 184 | 1097 | 0.2245277476 | 0.2243887540 |
| 185 | 1103 | 0.2243518169 | 0.2242140496 |
| 186 | 1109 | 0.2241771210 | 0.2240405627 |
| 187 | 1117 | 0.2239465673 | 0.2238111142 |
| 188 | 1123 | 0.2237746859 | 0.2236404092 |
| 189 | 1129 | 0.2236039887 | 0.2234708719 |
| 190 | 1151 | 0.2229916331 | 0.2228589972 |
| 191 | 1153 | 0.2229354535 | 0.2228041176 |
| 192 | 1163 | 0.2226619183 | 0.2225315407 |
| 193 | 1171 | 0.2224449867 | 0.2223156357 |
| 194 | 1181 | 0.2221768161 | 0.2220484009 |
| 195 | 1187 | 0.2220167854 | 0.2218894493 |
| 196 | 1193 | 0.2218577964 | 0.2217315248 |
| 197 | 1201 | 0.2216478323 | 0.2215225361 |
| 198 | 1213 | 0.2213368572 | 0.2212123766 |
| 199 | 1217 | 0.2212332446 | 0.2211098630 |
| 200 | 1223 | 0.2210792693 | 0.2209568994 |
| 201 | 1229 | 0.2209262660 | 0.2208048943 |
| 202 | 1231 | 0.2208746833 | 0.2207544371 |
| 203 | 1237 | 0.2207229689 | 0.2206036933 |
| 204 | 1249 | 0.2204235192 | 0.2203049964 |
| 205 | 1259 | 0.2201765851 | 0.2200588757 |
| 206 | 1277 | 0.2197392723 | 0.2196221029 |
| 207 | 1279 | 0.2196901674 | 0.2195740592 |
| 208 | 1283 | 0.2195933861 | 0.2194782596 |
| 209 | 1289 | 0.2194494963 | 0.2193352743 |
| 210 | 1291 | 0.2194010042 | 0.2192878017 |
| 211 | 1297 | 0.2192582684 | 0.2191459463 |
| 212 | 1301 | 0.2191632326 | 0.2190518416 |
| 213 | 1303 | 0.2191153390 | 0.2190049278 |
| 214 | 1307 | 0.2190208810 | 0.2189113760 |
| 215 | 1319 | 0.2187417521 | 0.2186329027 |
| 216 | 1321 | 0.2186947105 | 0.2185868054 |
| 217 | 1327 | 0.2185561371 | 0.2184490470 |
| 218 | 1361 | 0.2177902581 | 0.2176831719 |
| 219 | 1367 | 0.2176568353 | 0.2175505538 |
| 220 | 1373 | 0.2175241636 | 0.2174186770 |
| 221 | 1381 | 0.2173487375 | 0.2172439813 |
| 222 | 1399 | 0.2169599007 | 0.2168555965 |
| 223 | 1409 | 0.2167462119 | 0.2166425721 |
| 224 | 1423 | 0.2164506441 | 0.2163475570 |
| 225 | 1427 | 0.2163662251 | 0.2162639464 |
| 226 | 1429 | 0.2163236821 | 0.2162222530 |
| 227 | 1433 | 0.2162397304 | 0.2161390892 |
| 228 | 1439 | 0.2161148104 | 0.2160148969 |
| 229 | 1447 | 0.2159495744 | 0.2158503305 |
| 230 | 1451 | 0.2158669648 | 0.2157684816 |
| 231 | 1453 | 0.2158253479 | 0.2157276649 |
| 232 | 1459 | 0.2157026327 | 0.2156056434 |
| 233 | 1471 | 0.2154599545 | 0.2153635079 |
| 234 | 1481 | 0.2152595043 | 0.2151636431 |
| 235 | 1483 | 0.2152189718 | 0.2151238764 |
| 236 | 1487 | 0.2151389327 | 0.2150445476 |
| 237 | 1489 | 0.2150986227 | 0.2150049850 |
| 238 | 1493 | 0.2150190072 | 0.2149260626 |
| 239 | 1499 | 0.2149004870 | 0.2148081826 |
| 240 | 1511 | 0.2146660191 | 0.2145742151 |
| $m$ | $p_m$ | $t_m$ | $t_0$ |